\documentclass{amsart} 
\usepackage{latexsym} 
\usepackage{amssymb} 
\usepackage{amsmath}
\usepackage{enumerate}
\usepackage{hyperref}  

\begin{document}


\newtheorem{theorem}{Theorem} 
\newtheorem{problem}{Problem} 
\newtheorem{definition}{Definition} 
\newtheorem{lemma}{Lemma} 
\newtheorem{proposition}{Proposition} 
\newtheorem{corollary}{Corollary} 
\newtheorem{example}{Example} 
\newtheorem{conjecture}{Conjecture} 
\newtheorem{algorithm}{Algorithm} 
\newtheorem{exercise}{Exercise} 
\newtheorem{remarkk}{Remark} 
 
\newcommand{\be}{\begin{equation}} 
\newcommand{\ee}{\end{equation}} 
\newcommand{\bea}{\begin{eqnarray}} 
\newcommand{\eea}{\end{eqnarray}} 

\newcommand{\eeq}{\end{equation}} 

\newcommand{\eeqn}{\end{eqnarray}} 
\newcommand{\beaa}{\begin{eqnarray*}} 
\newcommand{\eeaa}{\end{eqnarray*}} 

\newcommand{\lip}{\langle} 
\newcommand{\rip}{\rangle}

\newcommand{\uu}{\underline} 
\newcommand{\oo}{\overline} 
\newcommand{\La}{\Lambda} 
\newcommand{\la}{\lambda} 
\newcommand{\eps}{\varepsilon} 
\newcommand{\om}{\omega} 
\newcommand{\Om}{\Omega} 
\newcommand{\ga}{\gamma} 
\newcommand{\rrr}{{\Bigr )}} 
\newcommand{\qqq}{{\Bigl\|}} 
 
\newcommand{\dint}{\displaystyle\int} 
\newcommand{\dsum}{\displaystyle\sum} 
\newcommand{\dfr}{\displaystyle\frac} 
\newcommand{\bige}{\mbox{\Large\it e}} 
\newcommand{\integers}{{\Bbb Z}} 
\newcommand{\rationals}{{\Bbb Q}} 
\newcommand{\reals}{{\rm I\!R}} 
\newcommand{\realsd}{\reals^d} 
\newcommand{\realsn}{\reals^n} 
\newcommand{\NN}{{\rm I\!N}} 
\newcommand{\DD}{{\rm I\!D}} 
\newcommand{\degree}{{\scriptscriptstyle \circ }} 
\newcommand{\dfn}{\stackrel{\triangle}{=}} 
\def\complex{\mathop{\raise .45ex\hbox{${\bf\scriptstyle{|}}$} 
     \kern -0.40em {\rm \textstyle{C}}}\nolimits} 
\def\hilbert{\mathop{\raise .21ex\hbox{$\bigcirc$}}\kern -1.005em {\rm\textstyle{H}}} 
\newcommand{\RAISE}{{\:\raisebox{.6ex}{$\scriptstyle{>}$}\raisebox{-.3ex} 
           {$\scriptstyle{\!\!\!\!\!<}\:$}}} 
 
\newcommand{\hh}{{\:\raisebox{1.8ex}{$\scriptstyle{\degree}$}\raisebox{.0ex} 
           {$\textstyle{\!\!\!\! H}$}}} 

\newcommand{\OO}{\won} 
\newcommand{\calA}{{\mathcal A}} 
\newcommand{\calB}{{\cal B}} 
\newcommand{\calC}{{\cal C}} 
\newcommand{\calD}{{\cal D}} 
\newcommand{\calE}{{\cal E}} 
\newcommand{\calF}{{\mathcal F}} 
\newcommand{\calG}{{\cal G}} 
\newcommand{\calH}{{\cal H}} 
\newcommand{\calK}{{\cal K}} 
\newcommand{\calL}{{\mathcal L}} 
\newcommand{\calM}{{\mathcal M}} 
\newcommand{\calO}{{\cal O}} 
\newcommand{\calP}{{\cal P}} 
\newcommand{\calU}{{\mathcal U}} 
\newcommand{\calX}{{\cal X}} 
\newcommand{\calXX}{{\cal X\mbox{\raisebox{.3ex}{$\!\!\!\!\!-$}}}} 
\newcommand{\calXXX}{{\cal X\!\!\!\!\!-}} 
\newcommand{\gi}{{\raisebox{.0ex}{$\scriptscriptstyle{\cal X}$} 
\raisebox{.1ex} {$\scriptstyle{\!\!\!\!-}\:$}}} 
\newcommand{\intsim}{\int_0^1\!\!\!\!\!\!\!\!\!\sim} 
\newcommand{\intsimt}{\int_0^t\!\!\!\!\!\!\!\!\!\sim} 
\newcommand{\pp}{{\partial}} 
\newcommand{\al}{{\alpha}} 
\newcommand{\sB}{{\cal B}} 
\newcommand{\sL}{{\cal L}} 
\newcommand{\sF}{{\cal F}} 
\newcommand{\sE}{{\cal E}} 
\newcommand{\sX}{{\cal X}} 
\newcommand{\R}{{\rm I\!R}} 
\renewcommand{\L}{{\rm I\!L}} 
\newcommand{\vp}{\varphi} 
\newcommand{\N}{{\rm I\!N}} 
\def\ooo{\lip} 
\def\ccc{\rip} 
\newcommand{\ot}{\hat\otimes} 
\newcommand{\rP}{{\Bbb P}} 
\newcommand{\bfcdot}{{\mbox{\boldmath$\cdot$}}} 
 
\renewcommand{\varrho}{{\ell}} 
\newcommand{\dett}{{\textstyle{\det_2}}} 
\newcommand{\sign}{{\mbox{\rm sign}}} 
\newcommand{\TE}{{\rm TE}} 
\newcommand{\TA}{{\rm TA}} 
\newcommand{\E}{{\rm E\, }} 
\newcommand{\won}{{\mbox{\bf 1}}} 
\newcommand{\Lebn}{{\rm Leb}_n} 
\newcommand{\Prob}{{\rm Prob\, }} 
\newcommand{\sinc}{{\rm sinc\, }} 
\newcommand{\ctg}{{\rm ctg\, }} 
\newcommand{\loc}{{\rm loc}} 
\newcommand{\trace}{{\, \, \rm trace\, \, }} 
\newcommand{\Dom}{{\rm Dom}} 
\newcommand{\ifff}{\mbox{\ if and only if\ }} 
\newcommand{\nproof}{\noindent {\bf Proof:\ }} 
\newcommand{\nproofYWN}{\noindent {\bf Proof of Theorem~\ref{YWN}:\ }} 
\newcommand{\remark}{\noindent {\bf Remark:\ }} 
\newcommand{\remarks}{\noindent {\bf Remarks:\ }} 
\newcommand{\note}{\noindent {\bf Note:\ }}

\newcommand{\boldx}{{\bf x}} 
\newcommand{\boldX}{{\bf X}} 
\newcommand{\boldy}{{\bf y}} 
\newcommand{\boldR}{{\bf R}} 
\newcommand{\uux}{\uu{x}} 
\newcommand{\uuY}{\uu{Y}} 
 
\newcommand{\limn}{\lim_{n \rightarrow \infty}} 
\newcommand{\limN}{\lim_{N \rightarrow \infty}} 
\newcommand{\limr}{\lim_{r \rightarrow \infty}} 
\newcommand{\limd}{\lim_{\delta \rightarrow \infty}} 
\newcommand{\limM}{\lim_{M \rightarrow \infty}} 
\newcommand{\limsupn}{\limsup_{n \rightarrow \infty}} 
 
\newcommand{\ra}{ \rightarrow } 

 \newcommand{\mlim}{\lim_{m \rightarrow \infty}}  
 \newcommand{\limm}{\lim_{m \rightarrow \infty}}  
 \newcommand{\nlim}{\lim_{n \rightarrow \infty}} 
 
 
 
 
 
 
 
\newcommand{\one}{\frac{1}{n}\:} 
\newcommand{\half}{\frac{1}{2}\:} 
 
\def\le{\leq} 
\def\ge{\geq} 
\def\lt{<} 
\def\gt{>} 
 
\def\squarebox#1{\hbox to #1{\hfill\vbox to #1{\vfill}}} 
\newcommand{\nqed}{\hspace*{\fill} 
           \vbox{\hrule\hbox{\vrule\squarebox{.667em}\vrule}\hrule}\bigskip} 

\title[Causal transference plans and their Monge-Kantorovich problems]{Causal transference plans and their Monge-Kantorovich problems}

\author{R\'{e}mi Lassalle}

\maketitle

\noindent 
{\bf Abstract : }{\small{This paper investigates causal optimal  transportation problems, in the framework of two Polish spaces, both endowed with filtrations. Specific concretizations yield primal problems equivalent to several classical problems of stochastic control, and of stochastic calculus ;   trivial filtrations yield usual problems of optimal transport.   Within this framework, primal attainments and dual formulations are obtained, under standard hypothesis,  for the related variational problems.  These problems are intrinsically related to martingales. Finally, we investigate applications to  stochastic frameworks. A straightforward equivalence between specific causal optimization problems, and problems of stochastic control, is obtained. Solutions to a class of stochastic differential equations are characterized, as optimum to specific causal Monge-Kantorovich problems ; the existence of a unique strong solution is related to corresponding Monge problems.}}
\\ 

\vspace{0.5cm}

 
\noindent 

\section*{{}} Over the last decade, connections between optimal transport and stochastic calculus have received  contributions of several origins  (among many, see  \cite{Soner},  \cite{leo},   \cite{Mika3}, \cite{Touzi}), with applications to fields such as \textit{financial mathematics} and \textit{stochastic mechanics}.  Here we establish, and we investigate, an extension of classical optimal transport, which encompasses applications  to \textit{stochastic differential equations}, and to \textit{stochastic control}.  

 Given two Polish spaces $E$ and $S$, optimal transport models transformations of a Borel probability  $\eta\in \mathcal{P}_E$,  to a probability $\nu \in \mathcal{P}_S$, by Borel probabilities on $E\times S$, whose first (resp. second) marginal is $\eta$ (resp. $\nu$). The latter set, denoted by $\Pi(\eta,\nu)$, is called the set of \textit{transference plans} (or of coupling plans), from  $\eta$ to $\nu$. To take into account the \textit{arrow of time}, in this paper, we endow $E$ (resp. $S$) with any \textit{filtration $(\mathcal{B}_t(E))_{t\in I }$ (resp. $(\mathcal{B}_t(S))_{t\in I}$), of its Borel sigma-field, indexed by a same totally ordered set $I$}.  As stated accurately below, within this framework, any transference plan $\gamma \in \Pi(\eta,\nu)$ induces canonically a filtration $(\mathcal{G}_t(\gamma))$ on $E$. We say that $\gamma$ is \textit{causal} if $$(\mathcal{G}_t(\gamma)) \subset( \mathcal{B}_t(E)^\eta), $$ the name being inherited from adapted, also called causal, processes ; subsequently, $\Pi_c(\eta,\nu)$ denotes the set of causal transference plans from $\eta$ to $\nu$. Roughly speaking, in most applications, a \textit{causal transference plan is a  transference plan, such that, at any time $t\in I$, the proportion of mass which is transported to any subset  $A\in \mathcal{B}_t(S)$, of the target space $S$, can be computed from the information available, at time $t$, on the initial space $E$} ; this information being modeled by the $\eta-$completion of $\mathcal{B}_t(E)$.  In this paper we consider \textit{optimal transportation problems} under this constraint. 

The structure of the paper is divided in two parts. The first part (Section~\ref{1}-Section~\ref{3}) investigates \textit{causal counterparts to  Monge-Kantorovich problems}, and related \textit{optimization} problems, in the \textit{analytic framework} stated above. Section~\ref{1} introduces the notation, and definitions, used in the whole paper.  Topological properties of \textit{causal transference plans} (Definition~\ref{defmain})  are stated in Theorem~\ref{theoremcompact}.   Without further assumptions, neither on marginals, nor  on filtrations, given a non-negative l.s.c. (lower semi-continuous) cost map $c : E\times S \to \mathbb{R}\cup\{+\infty\},$  we obtain the \textit{primal attainment} (Corollary~\ref{corolltp}), for the \textit{causal Monge-Kantorovich problem} $$P_{\eta,\nu} = \inf\left(\left\{ \int_{E\times S} c(x,y) d\gamma(x,y) \middle| \gamma\in \Pi_c(\eta,\nu)\right\} \right),$$  and the precise \textit{dual problem} (Theorem~\ref{thmduality}),   \begin{equation} \label{dualintro} P_{\eta,\nu} = \sup\left(\left\{ \int_S g(y) d\nu(y) \middle|  (g ,h) \in C_b(S)\times\mathcal{H}_{\eta,\nu} \ : \  h(x,y) + g(y) \leq c(x,y), \  \forall (x,y) \in E\times S\right\}\right), \end{equation}   $\mathcal{H}_{\eta,\nu}$ denoting a  set of maps which \textit{expresses} the causal \textit{constraint}, though a \textit{penalization} of the cost ;  this set is naturally related to \textit{martingales}.  The price to pay, for the generality of~(\ref{dualintro}), is Lemma~\ref{lemmacostdual2}, which extends the usual Kantorovich duality theorem, to a specific set of measurable cost maps, which are not lower-semicontinuous.

As a common feature, in most applications we encountered, the \textit{filtration on $S$ } satisfies further properties, emphasized explicitly in the proof of the celebrated Yamada-Watanabe criterion, on stochastic differential equations (see  Lemma 1.1. p.165 of \cite{I-W}) ; the  latter motivated this work. Namely, it is \textit{induced by a family $(\rho_t)_{t\in[0,1]}$ of continuous maps},  such that $\rho_t : S\to S$ is continuous, and $\mathcal{B}_t(S)= \rho_{t}^{-1}(\mathcal{B}(S))$, for all $t\in I$. It further satisfies the consistency condition  $$\rho_s\circ\rho_t = \rho_{s\wedge t},$$   for all  $s,t\in I$, $s\wedge t$ denoting the minimum of $s$ and $t$, and $\circ$ denoting the pullback of maps.  

To obtain compact statements,  we call \textit{regular}  (Definition~\ref{definitionconsist} of Section~\ref{3}) such filtrations ; it is  of topological origin. This involves further topological properties, inherited from those of the underlying space $S$,  for several sets of causal transference plans, with applications to  \textit{optimization problems} ; the latter are  investigated in Section~\ref{3}.  Lemma~\ref{consistentreg} is the key result of this section. It provides a characterization of {causal transference plans} within this assumption. Corollary~\ref{thmduality1} and Corollary~\ref{ToSEP}, which state \textit{primal attainment}  and \textit{dual problems} to several \textit{causal optimization} problems,  are obtained from Theorem~\ref{thmduality} ; their proofs emphasize the importance to have continuous and bounded maps in~(\ref{dualintro}), which is the main technical issue.  

This section enlightens the origin of applications to {stochastic control}, and to {stochastic calculus} investigated subsequently. Indeed,  Lemma~\ref{consistentreg} states that, within this further regularity assumption on $(\mathcal{B}_t(S))$, a \textit{causal transference plan} is a transference plans, which satisfies the constraint on joint laws of pairs of random elements, introduced by Yamada and Watanabe,  in the proof of the above mentioned result (see p.164-166 of \cite{I-W}). Equivalently, in this case, it is a \textit{transference plan, which preserves martingales}, in the precise sense that, the filtered probability space $(E\times S, (\mathcal{B}(E)\otimes \mathcal{B}(S))^\gamma, ((\mathcal{B}_t(E)\otimes \mathcal{B}_t(S))^\gamma),\gamma)$ is an \textit{extension of the filtered probability space} $(E ,\mathcal{B}(E)^\eta, (\mathcal{B}_t(E)^\eta), \eta)$ ; see Definition 7.1 of \cite{I-W}.  That is,  for any $(\mathcal{B}_t(E)^\eta)-$martingale $(M_t)_{t\in I}$, of the particular form $M_t= E_\eta[f| \mathcal{B}_t(E)^\eta]$, $f\in \mathcal{L}^\infty(\eta)$, on the complete probability space $(E ,\mathcal{B}(E)^\eta,\eta)$, the stochastic process $(M_t\circ \pi)_{t\in I}$  coincides with the  $(\mathcal{B}_t(E)\otimes \mathcal{B}_t(S))^\gamma-$martingale  $(E_\gamma[f\circ \pi | (\mathcal{B}_t(E)\otimes \mathcal{B}_t(S))^\gamma])$, on the complete probability space $(E\times S,\mathcal{B}(E\times S)^\gamma, \gamma)$; $\pi$ (resp. $\circ$) denoting the canonical projection $\pi : E\times S \to E$ (resp. the pullback of maps). Thus, \textit{causal optimal transportation problems} are \textit{intrinsically related to martingales}, and roughly speaking, within this further specific assumption, the originality of the present paper is to consider optimal mass transportation problems, over transference plans satisfying this Yamada-Watanabe constraint, or equivalently, preserving those martingales ; this seams to be new.

 The second part of the paper (Section~\ref{4}-Section~\ref{6}) investigates equivalences of causal optimization problems, to problems of stochastic calculus, and of stochastic control.  Section~\ref{4} is essentially a reformulation of part of the proof of the Yamada-Watanabe criterion, and of the related results of \cite{Jacod2}, within the analytic framework of the first part of the paper. Assuming the \textit{weak existence} and \textit{weak uniqueness} of solutions, Proposition~\ref{propja1} states that, taking suitable spaces and filtrations, given $X$, a weak solution to  \begin{equation} \label{SDE}dX_t = \sigma_t(X) dB_t + b_t(X)dt ; X_0=x, \end{equation}  on a complete stochastic basis $(\Omega, \mathcal{A}, (\mathcal{A}_t)_{t\in[0,1]}, \mathcal{P})$, and $\tau$ a $(\mathcal{A}_t)-$stopping-time (bounded by $1$), resp. $(Y_t)$ an $(\mathcal{A}_t)-$adapted continuous process, the joint law of the pair $(X,\tau)$ (resp. $(X,Y)$) is a causal transference plan, of the unique law $\eta$ of solutions to~(\ref{SDE}). Conversely, Proposition~\ref{propja2} states that, within the same framework, any causal transference plan of $\eta$, to any Borel probability on $[0,1]$ (resp. on $C([0,1], \mathbb{R}^d)$), is the joint law of such pairs.  We emphasize here that this equivalence fails, in general, when $(\mathcal{A}_t)$ is not right-continuous, which motivated our framework in subsequent sections. Section~\ref{5} (resp. Section~\ref{6}) investigates applications of \textit{causal optimization} problems (resp. \textit{causal Monge-Kantorovich problems}) to  \textit{stochastic control} (resp. to  \textit{stochastic differential equations}). By Section~\ref{4}, the equivalence is obtained in a systematic way : the processes, or stopping times, of interest in stochastic frameworks, are represented by the canonical projections of $E\times S$, which is straightforward. 
 
 Under the same assumptions on~(\ref{SDE}) as above,  Section~\ref{5} investigates the associated \textit{optimal stopping problems},  of the form \begin{equation}  \label{osint} \inf E_{\mathcal{P}}\left[c_a([X]_\tau, \tau)\right] ,\end{equation}  where the optimization is performed on all pairs $(X,\tau)$, defined on any complete stochastic basis $(\Omega, \mathcal{A}, ( \mathcal{A}_t)_{t\in[0,1]}, \mathcal{P})$, such that $\tau$ (resp. $X$) is an $(\mathcal{A}_t)-${stopping time} bounded by $1$ (resp. is a solution to~(\ref{SDE})), and where $[X]_\tau$ denotes the process $X$ stopped at $\tau$. Under weak conditions, we obtain the equivalence of~(\ref{osint}),  with the \textit{causal optimization} problem $$\inf\left(\left\{ \int_{W\times [0,1]} c(x,y) d\gamma(x,y) \middle| \gamma \in \mathcal{P}_c^{[0,1]}(\eta) \right\}\right),$$ where $\mathcal{P}_c^{[0,1]}(\eta) := \cup_{\nu\in \mathcal{P}_{[0,1]}}\Pi_c(\eta, \nu)$ denotes causal transference plans, from the law $\eta$ of solutions of~(\ref{SDE}), to any Borel probability on $[0,1]$ (the latter being endowed with a suitable filtration), and where $$c(\omega,t):= c_a(\omega_{.\wedge t},t),$$ for all $(t,\omega)\in [0,1] \times C([0,1], \mathbb{R}^d)$.This provides the \textit{primal attainment} for~(\ref{osint}), under very weak assumptions.   We relate this to the concomittent paper \cite{BHC}.  The latter, investigates an analytic approach to \textit{optimal SEP} (\textit{Skorokhod embedding problems}) ;  in their particular framework, causal transference plans appear, through objects they call $RST$ and $PRST$, and the equivalence is achieved by an auxiliary stopping time ; they obtain highly impressive results, through a sharp investigation of the geometry of their optimum. Here, we show that the same straightforward representation by projections, which we encountered in optimal stopping problems, and which we will meet again in applications to stochastic differential equations, yields a straightforward alternative \textit{representation} to \textit{extended optimal SEP}, which improves \cite{BHC} ; our representation is straightforward, and in view of applications to the primal attainment, we are more general. We emphasize here that \cite{BHC} motivated us to add general dual formula to this paper. Thus, the related sections are indebted to \cite{BHC}, and must be compared to it. As far as dual formulas in Polish frameworks are concerned, in a different framework, which yields more general applications, we follow a similar strategy ; in applications, the alternative approach we introduce here, which is natural,  involves specific technical difficulties described below. In particular the whole generality of~(\ref{dualintro}), which allows applications to non-Markovian frameworks, requires the specific extension, of the classical Kantorovich duality, by Lemma~\ref{lemmacostdual2}. Within the framework of Section~\ref{4},  Corollary~\ref{ToSEP} can be used to obtain the particular martingale dual formulas of \cite{BHC},  for generally non-Markovian solutions of (\ref{SDE}), under the assumptions of Section~\ref{4} ; this is not investigated in this paper.
 
  In Section~\ref{6}, under technical conditions, setting $W:= C([0,1],\mathbb{R}^d)$, and taking suitable filtrations,  we characterize (Theorem~\ref{SDEMK}) the joint law of pairs of processes $(X,B)$, \textit{weak solutions} to \textit{stochastic differential equations} (SDE) of the form  \begin{equation} \label{mpel4321} dX_t =dB_t +b_t(X) dt ; X_0=0, \end{equation}  as optimum of  \textit{causal Monge-Kantorovich problems} of the form $$\inf\left(\left\{ \int_{W\times W} |x-y|_H^2 d\gamma(x,y) \middle| \gamma \in \Pi_c(\mu,\nu) \right\}\right),$$ where $\mu$ denotes the law of \textit{standard Brownian motions}, where $\nu$ denotes the law of solutions to~(\ref{mpel4321}), which is assumed to be unique, and where $|.|_H$ denotes the so-called \textit{Cameron-Martin norm}. The \textit{pathwise uniqueness} for~(\ref{mpel4321}) is then related to the corresponding \textit{causal Monge problem}. The proof relies on an optimal transport formulation, of a well known representation of the \textit{relative entropy} with respect to $\mu$ (Lemma~\ref{lemmaweak}), inherited from a celebrated result of F\"{o}llmer (\cite{F1}) ; Lemma~\ref{lemmaweak} completes several recent results, mentioned in Section~\ref{6}. Moreover,  the characterization follows from a straightforward representation of $B$ (resp. of $X$), by the respective projection on the product space $W\times W$ ; in full consistency with Section~\ref{5}.  By completing  the final version of this manuscript, we were acquainted of the related \cite{KURTZ}, of stochastic calculus ; the analytic optimal transportation framework of the present paper should be compared.  Excepted results on extensions of filtered probability spaces, the general functional analytic proofs, of the first part of our paper, still hold if we drop the condition that $I$ is totally ordered, and that the families $(\mathcal{B}_t(E))$ and $(\mathcal{B}_t(S))$ are increasing ; with no further condition. However in this case, on suitable spaces, one looses the property to preserve martingales, which we do not want, since it is the key in applications to stochastic frameworks we encountered. Finally, the geometry of the optimal transference plans we obtain, is not beyond the scope of this paper.

   \section{Preliminaries and notation}  \label{1}  
  \subsection{Framework} \textit{In the whole paper, $E$ and $S$ denote two Polish spaces, endowed with filtrations $(\mathcal{B}_t(E))_{t\in I}$ (resp. $(\mathcal{B}_t(S))_{t\in I}$) of their Borel sigma-fields  $\mathcal{B}(E)$ (resp. $(\mathcal{B}(S)$), indexed by a same totally ordered set $I$ ; that is, $\mathcal{B}_s(E) \subset \mathcal{B}_t(E) \subset \mathcal{B}(E)$ for all $s\leq t$, and similarly on $S$. Further assumptions on spaces, on filtrations, and on $I$, are stated explicitly, when it is required.} 

Given a Polish space $Z$, $\mathcal{P}_Z$ denotes the set of Borel probabilities on $Z$. We systematically endow it with the so-called topology of weak convergence in measure (see \cite{STROOCK}), inherited from the weak$-\star$ topology, by identifying $\mathcal{P}_Z$ to a subset of the dual of $C_b(Z)$ ; the latter is endowed with the topology of uniform convergence.  Given a probability space $(\Omega, \mathcal{A}, \mathcal{P})$, $f_\star \mathcal{P}\in \mathcal{P}_Z$ denotes the direct image (pushforward), of a probability $\mathcal{P}$, by a measurable map $f: \Omega \to Z$.
To handle specific problems investigated below, given a sigma-field $\mathcal{G}\subset \mathcal{B}(E)$, and  $\eta\in \mathcal{P}_E$, it is useful to consider its $\eta-$completion, denoted by $\mathcal{G}^\eta$ ; the unique extension of $\eta$ to $\mathcal{B}(E)^\eta$ is still denoted by $\eta$. Moreover, in applications, equivalences between variational problems, involve joint laws of pairs of maps, and the projections on $E\times S$. For this reason we need  the following notation. Given  $X : \Omega \to E$ (resp. $Y : \Omega \to S$) two measurable maps, defined on a same probability space $(\Omega, \mathcal{A}, \mathcal{P})$, $X\times Y$ denotes  the $\mathcal{A} / \mathcal{B}(E\times S)$ measurable map defined by  \begin{equation} \label{prodeflmprolt} X\times Y : \omega \in \Omega \to (X(\omega), Y(\omega))\in E\times S.\end{equation}  Thus, $(X\times Y)_\star \mathcal{P}$ is the joint law of the pair $(X,Y)$.   \textit{In the whole paper $\pi : ( x , y) \in E\times S \to x\in E$ (resp. $ \widetilde{\pi} : (x, y) \in E\times S \to y \in S$) denote the projections on the respective component of $E\times S$}.

\subsection{Transference plans and the related kernels}  Given $\eta\in \mathcal{P}_E, \nu\in \mathcal{P}_S$, define \begin{equation} \label{plannona}\Pi(\eta,\nu) = \left\{ \gamma \in \mathcal{P}_{E\times S} \middle| \pi_\star \gamma = \eta , \widetilde{\pi}_\star \gamma = \nu \right\},\end{equation}  the set of \textit{transference plans} from $\eta$ to $\nu$, and \begin{equation} \label{atp} \mathcal{P}^S(\eta) := \cup_{\nu\in \mathcal{P}_S} \Pi(\eta, \nu). \end{equation}   Recall that, by continuity of the projections, $\Pi(\eta,\nu)$ and $\mathcal{P}^S(\eta)$ are closed. Moreover, together with Prohorov's criterion on the Polish space $E\times S$, the inner regularity of marginals implies that $\Pi(\eta,\nu)$ is relatively compact, and therefore, compact. 

A straightforward way to see transference plans as models of \textit{mass transports}, is to introduce the \textit{desintegration kernel} of $\gamma\in \Pi(\eta,\nu)$, with respect to the projection on $E$. We state the definition ;  for further details see \cite{I-W} (Theorem 3.3. and p.164), or \cite{SV} (Chapter I), and  the references therein.  Since $E$ and $S$ are Polish spaces, for any $\gamma\in \Pi(\eta,\nu)$, there exists a map $$\Theta_\gamma :  (x, B) \in E \times \mathcal{B}(S) \to \Theta_\gamma^x(B) \in \mathbb{R},$$ which meets the following assumptions  \\ \textit{\begin{enumerate}[(i)] \item For any $x \in E$,  $\Theta_\gamma^x \in \mathcal{P}_S.$ 
\item For all $B\in \mathcal{B}(S)$, the map  $\phi_B : x \in E \to \Theta_\gamma^x(B) \in \mathbb{R}$ is $\mathcal{B}(E)/ \mathcal{B}(\mathbb{R})-$measurable.
\item For all $A\in \mathcal{B}(E),B \in \mathcal{B}(S)$, \begin{equation} \label{aloi} \gamma( A \times B) = \int_A \Theta_\gamma^x(B)d\eta(x). \end{equation} 
\end{enumerate}  } It is further said to be \textit{unique}, in the sense that, if $\widetilde{\Theta}$ is a kernel which satisfies the same conditions, then outside a $\eta-$negligible set,  $\widetilde{\Theta}$ necessarily coincides with $\Theta_\gamma,$ 
as $\mathcal{P}_S-$valued maps. Actually,\begin{equation} \label{Tkernel} \Theta^x_\gamma(B) = \gamma( \widetilde{\pi} \in B |  \pi =x ) \ \eta-a.s. , \end{equation} so that $\Theta_\gamma^x(B)$ is interpreted, as the proportion of mass located at $x\in E$, which is transported to $B\in \mathcal{B}(S)$. Subsequently, for convenience of notation,  \begin{equation} \label{kernot} \gamma = \int_E d\eta(x) \delta^{Dirac}_x\otimes \Theta_\gamma^x\end{equation} denotes that $\Theta_\gamma$ is a \textit{kernel} associated to a transference plan $\gamma$, as above.

\subsection{Deterministic transference plans, and morphisms of probability spaces} \textit{Deterministic transference plans} naturally model systems, which answer deterministically to a, possibly random, input ; such plans can be seen as a transport, without splitting the mass from a given point.  By following \cite{MAL2}, given $(\Omega, \mathcal{A}, \mathcal{P})$, a complete probability space, $M_{\mathcal{P}}((\Omega, \mathcal{A}), (S,\mathcal{B}(S)))$ (or, when there are no ambiguity on the underlying space, $L^0(\mathcal{P}, S)$), denotes the set of $\mathcal{P}-$equivalence classes of maps, which is obtained by identifying the $\mathcal{A} / \mathcal{B}(S)-$ measurable maps $f : \Omega \to S$, which coincide outside $\mathcal{P}-$negligible sets.   Taking $(\Omega, \mathcal{A}, \mathcal{P})$ to be the completed space $(E, \mathcal{B}(E)^\eta, \eta)$, we systematically endow  $L^0(\eta,S)$ with its usual metric topology (see \cite{MAL2}) ; the latter induces the convergence in probability.  Thus, we refer to the latter, as the \textit{topology of convergence in probability}. Given $\nu\in \mathcal{P}_S,$  define $$\mathcal{R}(\eta, \nu):=\left\{ U \in L^0(\eta,S) \middle| U_\star \eta = \nu \right\}.$$  By following \cite{MAL2}, we call elements of $\mathcal{R}(\eta,\nu)$ \textit{morphisms of probability spaces}. Given $U\in L^0(\eta,S)$, $(I_E\times U)_\star \eta$ denotes $(I_E\times f_U)_\star\eta$, for any (and then all) $\mathcal{B}(E)^\eta/ \mathcal{B}(S)-$measurable map $f_U: E\to S$, whose $\eta-$equivalence class of maps is $U$, $I_E: x\in E\to x\in E$ denoting the identity map on $E$.  The \textit{continuous} injection  \begin{equation} \label{jdef} j : U \in L^0(\eta,S)  \to (I_E\times U)_\star \eta \in \mathcal{P}_{E\times S} \end{equation}  yields the \textit{embeddings} $L^0(\eta,S )  \hookrightarrow_j \mathcal{P}^S(\eta),$  and  $\mathcal{R}(\eta,\nu) \hookrightarrow_j \Pi(\eta,\nu),$ of \textit{morphisms of probability spaces} in \textit{transference plans}.  The set of \textit{deterministic transference plans}, from $\eta$ to $\nu$, is the range $j(\mathcal{R}(\eta,\nu))$ of $j|_{\mathcal{R}(\eta,\nu)}$, the restriction of $j$ to $\mathcal{R}(\eta,\nu)$ ; below, this notation yields compact statements, and proofs. Whence, for any $U\in \mathcal{R}(\eta,\nu)$,  $$(x,A) \in E\times \mathcal{B}(S) \to \delta^{Dirac}_{f_U(x)}(A):=1_A(f_U(x)) \in \mathbb{R},$$ is a \textit{kernel} associated to $\gamma:=j(U)$ by~(\ref{kernot}), for any $\mathcal{B}(E)^\eta/\mathcal{B}(S)-$measurable map $f_U: E\to S$, whose $\eta-$equivalence class is $U$. 

In \textit{stochastic analysis}, transformations of measures are usually achieved, on suitable spaces, as {direct images of laws of stochastic processes by morphisms of probability spaces,  which satisfy some further constraints}. In particular, those  induced by \textit{adapted processes}  are crucial. Causal transference plans, investigated below, can be seen as an abstraction of the related transference plans, by further allowing the mass to be splitted ; it will be clear in Section~\ref{6}.

\section{{Causal transference plans, and their optimal transportation problems}} \label{2} As stated accurately by Definition~\ref{defmain}, within the \textit{Polish framework} stated in Section~\ref{1}, once a \textit{transference plan} is given, it \textit{induces canonically a filtration} on the initial space $E$. The plan is \textit{causal} if its filtration is smaller than the completed reference filtration on $E$.  Definition~\ref{defmain} is natural, stemming from the functional analytic formulation (Lemma~\ref{convfprod}) of the Portmanteau theorem on product spaces ; it yields compact sets  (Theorem~\ref{theoremcompact}). Whence, under usual assumptions, the \textit{primal attainment} is obtained (Corollary~\ref{corolltp}), for causal Monge-Kantorovich problems. The natural idea, to handle dual problems, is to penalize cost maps, in order to apply the classical \textit{Kantorovich duality} theorem, through a \textit{min-max} theorem. A specific set of maps (Definition~\ref{heta}) provides a suitable penalization of the cost (Lemma~\ref{pmropo3}). However it is not l.s.c., so that the usual \textit{Kantorovich duality} doesn't apply.  General extensions of the latter require further conditions on marginals, which we do not want. We rather extend the classical Kantorovich duality to specific cost maps involved in the penalization (Lemma~\ref{lemmacostdual2}), by \textit{Lusin's theorem}. This yields the desired duality formula (Theorem~\ref{thmduality}), for arbitrary filtrations of the Borel sigma-fields, under the standard hypothesis of optimal transport ; none further condition is required, neither on filtrations, nor on on marginals.

\subsection{Definition, topological properties}  \begin{definition}
\label{defmain}
Given $\eta\in\mathcal{P}_E$, and $\nu \in \mathcal{P}_S$,  let $\gamma \in \Pi(\eta, \nu)$. For all  $t\in I$,  $\mathcal{G}_t(\gamma)$ denotes the $\eta-$completion, of the smallest sigma-field on $E$ such that, for any $C\in\mathcal{B}_t(S)$ of $\nu-$continuity (i.e. its boundary satisfies $\nu(\partial C)=0$), the map \begin{equation} \label{phicdef} \phi_C : x \in E \to \Theta_\gamma^x(C) \in \mathbb{R}\end{equation} is measurable, $\Theta_\gamma$  denoting any kernel, associated to $\gamma$ by~(\ref{kernot}).  We call $(\mathcal{G}_t(\gamma))_{t\in I}$, the filtration generated by $\gamma$ on $E$.  The set of causal transference plans (resp. of \textit{causals morphisms of filtered probability spaces}), from $\eta$ to $\nu$, denoted by $\Pi_c(\eta,\nu)$ (resp. by $\mathcal{R}_a(\eta,\nu)$), is defined by $$\Pi_c(\eta, \nu):=\left\{ \gamma \in \Pi(\eta,\nu) \middle|  \mathcal{G}_t(\gamma) \subset \mathcal{B}_t(E)^\eta \ ,  for \ all \ t \in I \right\},$$  resp. by $$\mathcal{R}_a(\eta,\nu):= \left\{ U \in \mathcal{R}(\eta,\nu) \middle|  j(U) \in \Pi_c(\eta,\nu)\right\},$$ $j$ denoting the map~(\ref{jdef}). 
\end{definition}
\begin{remarkk}
\begin{enumerate}[(i)]
\item Above, $(\mathcal{G}_t(\gamma))$ does not depend on the version of the kernel, associated to $\gamma \in \Pi(\eta,\nu)$ by~(\ref{kernot}). Whence the completion.
\item  From Definition~\ref{defmain}, $\gamma\in \Pi(\eta,\nu)$ is causal, if and only if, the map $\phi_C$ defined by~(\ref{phicdef}) is $\mathcal{B}_t(E)^\eta-$measurable, for any $C\in\mathcal{B}_t(S)$ of $\nu-$continuity, for any $t\in I$.
\item  If $Y :\Omega \to E$, and $Z : \Omega \to S$, are two  $\mathcal{A} / \mathcal{B}(E)-$ (resp.  $\mathcal{A} / \mathcal{B}(S)-$) measurable maps, defined on a same complete probability space $(\Omega, \mathcal{A}, \mathcal{P})$, such that $(Y\times Z)_\star\mathcal{P}  \in \Pi_c(Y_\star \mathcal{P},Z_\star \mathcal{P}),$ then $(Y, Z ,(\Omega, \mathcal{A}, \mathcal{P}))$ (or $(Y,Z)$ for short) will be called a causal coupling.
\end{enumerate}
\end{remarkk}
Recall that, given $\eta\in \mathcal{P}_E$,  $\mathcal{P}^S(\eta)$, defined by~(\ref{atp}), is the set of Borel probabilities on $E\times S$, whose first marginal is $\eta$. The following Lemma~\ref{convfprod} is a functional analytic formulation of the Pormanteau theorem on product spaces ; it yields the compactness of causal transference plans with given marginals.

 \begin{lemma}
  \label{convfprod}
  Given $\eta\in \mathcal{P}_E$, $\nu \in \mathcal{P}_S$, and $\gamma \in \Pi(\eta,\nu)$, a sequence $(\gamma_n)_{n\in \mathbb{N}}\subset \mathcal{P}^S(\eta)$ converges to $\gamma$, in the topology of weak convergence in measure, if and only if, for any Borel set $B\in \mathcal{B}(S)$ of $\nu-$continuity (i.e. $\nu(\partial B) = 0$),  the sequence $(\Theta_n(B))_{n\in \mathbb{N}}$ converges to $\Theta(B)$, in the weak topology $\sigma(L^1(\eta), L^\infty(\eta))$ of $L^1(\eta)$, for any (and then for all) kernel $\Theta$ (resp., for $n  \in \mathbb{N}$, $\Theta_n$) associated to $\gamma$ (resp.  to $\gamma_n$) by~(\ref{kernot}).
  \end{lemma}
 \nproof
Taking  a set $A\in \mathcal{B}(E)$ (resp. $B\in \mathcal{B}(S)$), of $\eta-$ (resp. $\nu-$) continuity, from the definition (see~(\ref{aloi})), the $\sigma(L^1(\eta), L^\infty(\eta))-
$convergence of $(\Theta_n(B))$ to $\Theta(B)$ yields $$\gamma(A\times B) = \limn \gamma_n(A\times B).$$  From the Portmanteau theorem on product spaces (for instance 
see Theorem 3.1 of \cite{BIL}), it implies the convergence of $(\gamma_n)$ to $\gamma$.  Conversely, denote by $\mathcal{C}_\eta$ the subset of the $g\in L^
\infty(\eta)$ of the form $\eta-a.s.$ \begin{equation} \label{gform} g= \sum_{i=1}^N \alpha_i 1_{A_i},\end{equation} for some $N\in \mathbb{N}$, where for all $i=1,...,N$, $\alpha_i\in \mathbb{R}$, and $A_i\in \mathcal{B}(E)$ is a set of $\eta-$continuity. $\mathcal{C}_\eta$ is dense in $L^\infty(\eta)$ for the (strong) $L^1(\eta)$ 
topology. Indeed, together with  Tietze's extension theorem, Lusin's theorem yields the density of $C_b(E)$ in $L^\infty(\eta)$, for the subspace topology, which is induced by the strong $L^1(\eta)-$topology. On the other hand, given $f\in C_b(E)$, by the dominated convergence theorem, a classical argument (for instance see the proof of Theorem 1.1.1 of \cite{SV}) ensures the existence of a sequence $(g_n)_{n\in \mathbb{N}} \subset \mathcal{C}_\eta$, which converges to $f$, in the strong $L^1(\eta)-$topology. Thus, given $\epsilon >0$, and $X\in L^\infty(\eta)$, we take $g\in \mathcal{C}_\eta$ of the form~(\ref{gform}), such that $$|X-g|_{L^1(\eta)} \leq \frac{\epsilon}{2}.$$ Given a set $B\in \mathcal{B}(S)$ of $\nu-$continuity,  we first obtain \begin{equation} \label{kerdef1MP} |E_\eta[X\Theta_n(B) ] -E_\eta[X 
\Theta(B)] | \leq \epsilon +  \sum_{i=1}^N |\alpha_i| |\gamma_n(A_i\times B) - \gamma(A_i\times B)|, \end{equation} where we used that {kernels} takes their values in $[0,1]$, together with the definition of $\Theta_n$, $\Theta$ (see~(\ref{aloi})), and $g$. Recall that $\eta$ (resp. $\nu$) are the first (resp. second) marginal of $\gamma$, and that $B$ (resp. for all $i\in [1,N]$, $A_i$) is a set of $\nu-$ (resp. $\eta-$) continuity. Therefore, together with~(\ref{kerdef1MP}), by the Portmanteau theorem on  product spaces, the weak convergence of $(\gamma_n)$ to $\gamma$ yields  \begin{equation} \label{majolem22} \limsup |E_\eta[X \Theta_n(B) - X \Theta(B)]| \leq \epsilon. \end{equation}  Since~(\ref{majolem22}) holds for any $\epsilon >0$, and for any $X\in L^\infty(\eta)$, the proof is complete.  \nqed 

Theorem~\ref{theoremcompact} states topological properties of the sets, introduced by Definition~\ref{defmain}.  \begin{theorem}
 \label{theoremcompact}
For any $\eta \in \mathcal{P}_E$, $\nu \in \mathcal{P}_S$, and any  filtration $(\mathcal{B}_t(E))_{t\in I}$ (resp. $(\mathcal{B}_t(S))_{t\in I}$) of the Borel sigma-field of $E$ (resp. of $S$),  indexed by a same totally ordered set $I$,  \begin{enumerate}[(i)] \item  $\Pi_c(\eta,\nu)$, the set of causal transference plans from $\eta$ to $\nu$,  is a not empty convex set, which is compact in $\mathcal{P}_{E\times S}$, for the topology of weak convergence in measure ; \item $\mathcal{R}_a(\eta, \nu)$, the set of causal morphisms of probability spaces from $\eta$ to $\nu$, is closed in $\mathcal{R}(\eta,\nu)$, endowed with the topology of convergence in probability (see Section~\ref{1}).    \end{enumerate}
\end{theorem}
  \nproof As it contains $\eta\otimes \nu$, $\Pi_c(\eta, \nu)$ is not empty, while the convexity trivially follows from the definition. $\mathcal{R}_a(\eta,\nu)=j^{-1}(\Pi_c(\eta,\nu))$, $j$ denoting the continuous map~(\ref{jdef}). Since any compact set of $\mathcal{P}_{E\times S}$ is closed, $(ii)$ follows from $(i)$.  Moreover, $\Pi(\eta,\nu)$ being compact, to obtain $(i)$, it is enough to prove that $\Pi_c(\eta,\nu)$ is closed. Thus, assume that $(\gamma_n)_{n\in \mathbb{N}} \subset \Pi_c(\eta, \nu)$, converges to some $\gamma \in \mathcal{P}_{E\times S}$. Since $\Pi(\eta,\nu)$ is closed, it contains $\gamma$ ; Given $t\in I$, Definition~\ref{defmain} yields \begin{equation} \label{thetanproof}  \Theta_n(B) = E_\eta\left[\Theta_n(B) \middle| \mathcal{B}_t(E)^\eta\right]   \ \eta-a.s., \end{equation}  for all $n\in \mathbb{N}$, and $B\in \mathcal{B}_t(S)$, which is a set of $\nu-$continuity ; $\Theta$ (resp.for all $n$, by $\Theta_n$), denoting any {kernel} associated to $\gamma$ (resp. to $\gamma_n$) by~(\ref{kernot}). On the other hand, Lemma~\ref{convfprod} implies the weak $L^1(\eta)$ convergence of $(\Theta_n(B))$ to $(\Theta(B))$, for all set $B\in \mathcal{B}(S)$ of $\nu-$continuity. Furthermore, as a linear map which is strongly continuous, the conditional expectation is weakly continuous in $L^1(\eta)$. Whence,  taking  the  $\sigma(L^1(\eta), L^\infty(\eta))-$limit in~(\ref{thetanproof}), it follows  that $\Theta(B)$ is $\mathcal{B}_t(E)^\eta-$measurable, which achieves the proof. \nqed
    \subsection{Primal attainment} Similarly to the usual \textit{Monge-Kantorovich problems},  together with the Portmanteau  theorem of \cite{STROOCK}, Theorem~\ref{theoremcompact} yields the  \textit{primal attainment} stated below.
    \begin{corollary}
 \label{corolltp} (Primal attainment)
  Given $$c:E\times S \to \mathbb{R} \cup \{+\infty\},$$  a non negative l.s.c. map, for any $\eta \in \mathcal{P}_E$, and $\nu \in \mathcal{P}_S$, there exists a  $\gamma \in \Pi_c(\eta,\nu)$, which attains the infimum  of \begin{equation} \label{eret} \inf\left(\left\{  \int_{E\times S} c(x,y) d\gamma(x,y) \middle| \gamma\in \Pi_c(\eta, \nu) \right\} \right). \end{equation}  \nqed
 \end{corollary}
We call~(\ref{eret}) a \textit{causal Monge-Kantorovich problem}, and we call \begin{equation} \label{causmm} P_{\eta,\nu}^{M}:=\inf\left(\left\{ \int_{E} c(x,U(x)) d\eta(x) \middle| U\in \mathcal{R}_a(\eta,\nu) \right\} \right), \end{equation} the \textit{causal Monge problem} associated to~(\ref{eret}).

\subsection{Dual problems} The following set of maps, closely related to \textit{martingales}, is naturally involved in dual formulations to \textit{causal Monge-Kantorovich problems}. \begin{definition}
\label{heta}
For convenience of notations, set $\mathcal{B}_{t_0}(S):= \{\emptyset, S\}$, and $\mathcal{B}_{t_0}(E):= \{\emptyset, E\}$, for some $t_0\notin I$. Given $\eta \in \mathcal{P}_E$, and $\nu\in \mathcal{P}_S$, we define $\mathcal{H}_{\eta,\nu}$ (resp. $\mathcal{H}_\eta$) to be the convex subset, of Borel measurable maps $h :E\times S \to \mathbb{R},$ of the specific form   \begin{equation} \label{hform} h=  \sum_{i=1}^N g_i \circ \pi 1_{A_i} \circ \widetilde{\pi}, \end{equation} where $N\in\mathbb{N}$,  and where, for some $(t_i)_{i=1,...,N} \subset I \cup\{t_0\} $, $(A_i)$, and $(g_i)$, meet the following assumptions, for all $i\in I\cup \{t_0\}$:
\begin{enumerate}[(i)]
\item  $A_i\in \mathcal{B}_{t_i}(S)$, and $\nu(\partial A_i) =0$ (resp. $A_i \in \mathcal{B}_{t_i}(S)$).
\item  $g_i: E\to \mathbb{R}$ is a bounded measurable map, such that  $\eta-a.s.$ \begin{equation} \label{refdauphf} X_i-E_{\eta}[X_i | \mathcal{B}_{t_i}(E)^\eta] =g_i, \end{equation} for some $X_i\in L^\infty(\eta)$.
\end{enumerate}
\end{definition}

\begin{lemma}
\label{pmropo3}
Given $\eta\in \mathcal{P}_E$, and $\nu \in \mathcal{P}_S$, define the map $$\mathcal{C}_{caus} : \gamma \in \Pi(\eta,\nu ) \to \mathcal{C}_{caus} (\gamma):=  \sup_{h\in \mathcal{H}_{\eta,\nu}} C_h(\gamma) \in \mathbb{R}\cup \{+\infty \}, $$ where for all $h\in \mathcal{H}_{\eta,\nu},$ $C_h$ denotes the map $$C_h : \gamma \in \Pi(\eta, \nu) \to C_h(\gamma) := E_\gamma[h] \in \mathbb{R}.$$ The following assertions hold :
\begin{enumerate}[(i)]
\item For all $h\in \mathcal{H}_{\eta,\nu}$, the map $C_h$ is continuous on $\Pi(\eta,\nu)$, for the topology of weak convergence in measure. 
\item For all $h\in \mathcal{H}_{\eta,\nu}$, and $\gamma \in \Pi_c(\eta,\nu)$, $E_\gamma[h] =0.$
\item  \begin{equation} \mathcal{C}_{caus}(\gamma) = \begin{cases}
       0 & \text{if } \ \gamma\in \Pi_c(\eta,\nu) \\
       + \infty & \text{otherwise}
       \end{cases} .\end{equation}
\end{enumerate}       
\end{lemma}
\nproof Taking  $h\in \mathcal{H}_{\eta,\nu}$, it is of the form~(\ref{hform}), where the $A_i$ are sets of $\nu-$continuity. Thus, for all $\gamma\in \Pi(\eta,\nu),$ we have  \begin{equation} \label{rplo} E_\gamma[h] = \sum_{i=1}^N E_\gamma[g_i \circ \pi  1_{A_i}\circ \widetilde{\pi} ]=  \sum_{i=1}
^NE_\eta\left[ g_i \Theta(A_i)\right], \end{equation} where $\Theta$ is a {kernel} associated to $\gamma$ by~(\ref{kernot}), and where, for all $i\in[1,N]$,~(\ref{refdauphf}) holds, for some $X_i\in L^\infty(\eta)$.  Given $(\gamma_n)_{n\in \mathbb{N}}\subset \Pi(\eta,\nu)$, a sequence of {transference plans}, which converges to some $\gamma\in \mathcal{P}_{E\times S}$,  since $\Pi(\eta,\nu)$ is closed (and compact), we obtain $\gamma\in \Pi(\eta,\nu)$. Moreover for all $i$, $A_i$ is a set of $\nu-$continuity. Whence, denoting by $\Theta_n$ the {kernel} associated to $\gamma_n$ by~(\ref{kernot}), for all $n\in \mathbb{N}$, Lemma~\ref{convfprod} implies the weak $L^1(\eta)-$convergence of $(\Theta_n(A_i))$ to $\Theta(A_i)$. Since $g_i \in L^\infty(\eta)$, by~(\ref{rplo}) it yields $$C_h(\gamma) = \sum_{i=1}^N E_\eta\left[g_i\Theta(A_i)\right] = \limn \sum_{i=1}^N   E_\eta[g_i \Theta_n(A_i) ]  =  \lim_n C_h(\gamma_n),$$ which yields $(i)$.   Taking $\gamma\in \Pi_c(\eta,\nu)$, and $\Theta$ an associated {kernel} by~(\ref{kernot}),  for all $h\in \mathcal
{H}_{\eta,\nu}$  of the form~(\ref{hform}),  by~(\ref{rplo}) we obtain $$E_\gamma[h]  = \sum_{i=1}^NE_\eta\left[ \left(X_i-E_{\eta}[X_i | \mathcal{B}_{t_i}(E)^\eta]\right) \Theta(A_i)\right] =0, $$ where we used that, since $\gamma$ is causal, and $A_i\in \mathcal{B}_{t_i}(S)$ is a set of $\nu-$continuity,  $\Theta(A_i)$ is $(\mathcal{B}_{t_i}(E)^\eta)$ measurable, for all $i\in I$.  To prove $(iii)$, first notice that, by $(ii)$, the map $\mathcal{C}_{caus}$ vanishes on $\Pi_c(\eta,\nu)$.  Conversely, taking $\gamma \in \Pi(\eta,\nu)$ such that $\mathcal{C}_{caus}(\gamma)=0$,  the linearity implies $E_\gamma[h] = 0,$ for all $h\in \mathcal{H}_{\eta,\nu}$. Given $s\in I$, and $B\in \mathcal{B}_s(S)$, a set of $\nu-$continuity, for all $A\in \mathcal{B}(E)$ we take $h\in \mathcal{H}_{\eta,\nu}$ such that $$h:= ( 1_A - E_{\eta}[1_A| \mathcal{B}_s(E)^\eta])\circ \pi 1_{B}\circ \widetilde{\pi}  \ \gamma-a.s..$$ Whence, we obtain  $$0 =  E_\gamma[h] =  E_\eta\left[\left(1_A - E_{\eta}\left[1_A\middle| \mathcal{B}_s(E)^\eta\right]\right) \Theta(B)\right] = E_\eta\left[1_A\left(\Theta(B) -E_\eta\left[\Theta(B)\middle| \mathcal{B}_s(E)^\eta\right]\right)\right].$$ The latter holds for all $A\in \mathcal{B}(E)$, so that $ \Theta(B)$ is $\mathcal{B}_t(E)^\eta-$measurable. This yields $\gamma \in \Pi_c(\eta,\nu)$ ; by linearity we obtain $(iii)$. \nqed
\begin{remarkk}
\label{remarkkFL2}
In particular, for all measurable  $c  : E\times S \to \mathbb{R} \cup\{+ \infty\}$,  Lemma~\ref{pmropo3} entails  $$P_{\eta,\nu} = \inf\left(\left\{E_\gamma[c] +  \mathcal{C}_{caus}(\gamma) \middle| \gamma \in \Pi(\eta,\nu) \right\}\right), $$ $P_{\eta,\nu}$ denoting the infimum of the {primal problem}~(\ref{eret}). \end{remarkk}
\begin{lemma}
\label{lemmacostdual2}
Given  a non-negative l.s.c. map $c :E\times S \to \mathbb{R} \cup\{+\infty\}$, let  $\eta \in \mathcal{P}_E$, and $\nu \in \mathcal{P}_S$. Further consider $h : E\times S \to \mathbb{R},$ a map which is assumed to be of the specific form $$h : (x,y) \in E\times S \to h(x,y) := \sum_{i=1}^N Z_i(x) 1_{A_i}(y)\in \mathbb{R},$$ for some  $N\in \mathbb{N}$,  where, $Z_i\in \mathcal{L}^\infty(\eta)$, and where $A_i\in \mathcal{B}(S)$ satisfies $\nu(\partial A_i) =0$, for all $i=1,...,N$. Defining $$D^{h}_{\eta,\nu}:= \sup\left(\left\{\int_E f(x) d\eta(x) + \int_S g(y) d\nu(y) \middle| (f,g)\in \mathcal{L}^\infty(\eta)\times C_b(S), f\circ \pi + g\circ \widetilde{\pi} \leq c+ h\right\}\right),$$ and $$P^{h}_{\eta,\nu}:= \inf\left(\left\{\int _{E\times S}(c(x,y) +h(x,y)) d\gamma(x,y) : \gamma \in \Pi(\eta,\nu) \right\}\right),$$  we have $P^{h}_{\eta,\nu} = D^{h}_{\eta,\nu}.$   \end{lemma}
\nproof Since $D^{h}_{\eta,\nu} \leq  P^{h}_{\eta,\nu}$  follows from the definition, it is enough to prove the converse inequality. Henceforth, take $\epsilon >0$, and for convenience of notation, define $B:=\max_{i=1,...,N}|Z_i|_{L^{\infty}(\eta)}$. Since  $Z_i\in \mathcal{L}^\infty(\eta)$, applying Lusin's theorem, together with Tietze's extension theorem,  we obtain the existence of $Z_i^{\epsilon}\in C_b(E)$, such that \begin{equation} \label{420tilde} |Z_i^{\epsilon}|_{L^\infty(\eta)} \leq B, \end{equation} for all $i=1,...,N$, and \begin{equation} \label{proofdp3} \max_{i=1,...,N} |Z_i - Z_i^{\epsilon}|_{L^1(\eta)} \leq \frac{\epsilon}{3 N}. \end{equation} On the other hand, $$1_{\overline{A}_i}(y) = \downarrow \limn \left(\frac{1}{1+d_S(y,\overline{A}_i)} \right)^{n},$$  for all $y\in S$, and $i=1,...,N$. Whence, the dominated convergence theorem ensures the existence of $g_i^\epsilon \in C_b(S)$, such that  \begin{equation} \label{421tilde} 1 \geq g_i^\epsilon \geq  1_{\overline{A_i}} \geq 1_{A_i}, \end{equation} and \begin{equation} \label{proofdp5} |g^\epsilon_i - 1_{\overline{A_i}}|_{L^1(\nu)} \leq \frac{\epsilon}{ 3 N B}, \end{equation} for all $i=1...,N$. Since $y\in S \to 1_{\overset{\circ}{A}_i}(y) \in \mathbb{R}$ is a non-negative l.s.c. map, together with the dominated convergence theorem,  Moreau-Yosida's approximations (see Lemma 4.4 of \cite{LEOCONV}) ensure the existence, of a Lipschitz continuous bounded map $G_i^\epsilon \in C_b(S)$, which satisfies  \begin{equation} \label{proofdp3} G^\epsilon_i \leq 1_{\overset{\circ}{A}_i} \leq 1_{A_i} \end{equation} and \begin{equation} \label{proofdp7} |G^\epsilon_i - 1_{\overset{\circ}{A}_i}|_{L^1(\nu)} \leq \frac{\epsilon}{ 3N B}, \end{equation} for all $i=1...,N$. Define $F_\epsilon := \sum_{i=1}^N |Z_i^\epsilon - Z_i|$, $G_\epsilon := B \sum_{i=1}^N (g_i^\epsilon - G_i^\epsilon)$, and $h_\epsilon : (x,y) \in E\times S \to h_\epsilon(x,y) := \sum_{i=1}^N Z_i^\epsilon(x)g_i^\epsilon(y)  \in \mathbb{R}.$ In particular, $(F_\epsilon, G_\epsilon,h_\epsilon) \in L^\infty(\eta)\times C_b(S) \times C_b(E\times S)$. From the definitions,~(\ref{420tilde}),~(\ref{421tilde}), and~(\ref{proofdp3}) yield  \begin{equation} \label{proofdpstar1} |h_\epsilon(x,y) - h(x,y)| \leq F_\epsilon(x) + G_\epsilon(y), \end{equation}  for all $(x,y)\in E\times S$. Since, $\nu(\partial A_i)=0$,  we obtain $1_{\overset{\circ}{A}_i}= 1_{\overline{A}_i} \ \nu-a.s.,$ for all $i=1,...,N$. Thus,~(\ref{proofdp5}), (\ref{proofdp3}), and (\ref{proofdp7})  yield \begin{equation}  \label{proofdpstar2} E_\eta[F_\epsilon]  + E_\nu[G_\epsilon]  \leq \epsilon. \end{equation} By~(\ref{proofdpstar1}) and~(\ref{proofdpstar2}), for any $\gamma\in \Pi(\eta,\nu)$, we obtain $$ P^{h}_{\eta,\nu} \leq E_\gamma[c+h] \leq E_\gamma[c+ h_\epsilon] + \epsilon, $$ so that \begin{equation} \label{proofdp7p} P^{h}_{\eta,\nu} \leq P^\epsilon + \epsilon, \end{equation} where $P^\epsilon := \inf_{\gamma\in \Pi(\eta,\nu)} E_\gamma[c+ h_\epsilon].$ By definition of $(g_i^\epsilon)$, $h^\epsilon$, and by~(\ref{420tilde}), the classical Kantorovich duality theorem (see Theorem 5.1 of \cite{Vil2}, or \cite{LEOCONV}) yields the existence of $(f,g)\in C_b(E)\times C_b(S)$ such that \begin{equation} \label{proofdp8} E_\eta[f] + E_\nu[g] \geq P^\epsilon -\epsilon, \end{equation} and  \begin{equation} \label{proofdp9} f(x) + g(y) \leq c(x,y) + h_\epsilon(x,y) ,\end{equation} for all $x\in E, y\in S$. Define $\widetilde{f}:= f- F_\epsilon \in L^\infty(\eta)$, and $\widetilde{g}:= g- G_\epsilon \in C_b(S)$. By~(\ref{proofdpstar1}) and~(\ref{proofdp9}) we obtain $$\widetilde{f}(x)+ \widetilde{g}(y) \leq  c(x,y) + h(x,y), $$  for all $(x,y)\in E\times S$. Whence,   $$D^{h}_{\eta,\nu}  \geq E_\eta[\widetilde{f}] + E_\nu[\widetilde{g}]\geq P^\epsilon - 2\epsilon  \geq P^{h}_{\eta,\nu} - 3\epsilon,$$ follows from~(\ref{proofdpstar2}),~(\ref{proofdp7p}), and~(\ref{proofdp8}), for all $\epsilon >0$ ; the proof is achieved. \nqed
\begin{theorem}
\label{thmduality}
Given  a non-negative l.s.c. map $c: E\times S \to \mathbb{R}\cup\{ +\infty\}$,  for $\eta \in \mathcal{P}_E$ and $\nu \in \mathcal{P}_S,$  we have  \begin{equation} \label{dualityabstrpm} P_{\eta,\nu} = {D}_{\eta,\nu}, \end{equation} where
 \begin{equation}
{D}_{\eta,\nu}:= \sup\left(\left\{\int_S g(y) d\nu(y) \middle| (g,h) \in C_b(S)\times \mathcal{H}_{\eta,\nu} :  \ g\circ \widetilde{\pi} + h \leq c \right\}\right), \end{equation} and  where $$P_{\eta,\nu}:= \inf\left(\left\{ \int_{E\times S} c(x,y) d\gamma(x,y) \middle| \gamma \in \Pi_c(\eta,\nu)\right\}\right).$$
\end{theorem}
\nproof By $(ii)$ of Lemma~\ref{pmropo3}, $D_{\eta,\nu}\leq P_{\eta,\nu}$. Since $\Pi_c(\eta,\nu)$ is compact (Theorem~\ref{theoremcompact}), taking Moreau-Yosida's approximation's (see Lemma 4.4 of \cite{LEOCONV}), it is enough to prove the converse inequality for  $c: E\times S \to \mathbb{R}$,  continuous and bounded.  Assuming that $c\in C_b(E\times S)$, Lemma~\ref{pmropo3} yields $$ P_{\eta,\nu} = \inf_{ \gamma \in \Pi(\eta,\nu)} \sup_{k\in  \mathcal{H}_{\eta,\nu}} F(\gamma, k)  \in \mathbb{R},$$ where $F : (\gamma, k)\in \Pi(\eta,\nu)  \times  \mathcal{H}_{\eta,\nu}  \to   E_\gamma[c + k] \in \mathbb{R}.$ For all  $\gamma\in \Pi(\eta,\nu)$ (resp. $k\in \mathcal{H}_{\eta,\nu}$), $F(\gamma,.)$ (resp. $F(.,k)$) is concave (resp. convex and l.s.c.), on the convex set $\mathcal{H}_{\eta,\nu}$ (resp. on the compact convex set $\Pi(\eta,\nu)$).  Indeed, by $(i)$ of Lemma~\ref{pmropo3} (resp. by the Portmanteau theorem of \cite{STROOCK}), the map $\gamma\in \Pi(\eta,\nu)  \to E_\gamma[k]\in \mathbb{R}$ (resp.  $\mathcal{C}: \gamma \in \mathcal{P}_{E\times S} \to \mathcal{C}(\gamma):=E_\gamma[c]\in \mathbb{R}$) is continuous, so that $F(.,k)$ is continuous. Therefore, by a classical result of convex analysis (for instance Theorem~2 of \cite{FAN}), we obtain \begin{equation} \label{mplf} P_{\eta,\nu} = \inf_{ \gamma \in \Pi(\eta,\nu)} \sup_{k\in  \mathcal{H}_{\eta,\nu}} F(\gamma,k) = \sup_{k
\in  \mathcal{H}_{\eta,\nu}}  \inf_{ \gamma \in \Pi(\eta,\nu)}  F(\gamma, k).\end{equation}  By  Lemma~\ref{lemmacostdual2}, given $\epsilon>0$,~(\ref{mplf}) ensures the existence of $k\in \mathcal{H}_{\eta,\nu}$, $f\in \mathcal{L}^\infty(\eta)$, $\widetilde{g} \in C_b(S)$ such that \begin{equation} \label{stardualthm} f(x) + \widetilde{g}(y) \leq c(x,y) + k(x,y), \end{equation} for all $(x,y)\in E\times S$, and  \begin{equation} \label{stardualthm2} E_\eta[f] + E_\nu[\widetilde{g}] \geq P_{\eta,\nu}-\epsilon. \end{equation} Define $h : (x,y)\in E\times S \to h(x,y) = f(x)-E_\eta[f]- k(x,y) \in \mathbb{R},$ and $g(y) = E_\eta[f] + \widetilde{g}(y).$ Thus, $h\in \mathcal{H}_{\eta,\nu}$, and $g\in C_b(S)$.  By~(\ref{stardualthm}) and~(\ref{stardualthm2}), it yields  $$D_{\eta,\nu} \geq E_\nu[g] =E_\eta[f] + E_\nu[\widetilde{g}] \geq P_{\eta,\nu} -\epsilon,$$ for all $\epsilon>0$, which proves the claim. \nqed 
\section{{Causal optimization problems, under regularity assumptions}} \label{3} Under further assumptions on $(\mathcal{B}_t(S))$, we consider several \textit{causal optimization} problems. To obtain compact statements, we call a \textit{regular filtration} (Definition~\ref{definitionconsist}), any filtration which meets some assumptions of topological origin, which are  satisfied in most situations encountered in applications. On the contrary, in view of applications, the filtration on the first space $(\mathcal{B}_t(E))$ is still assumed to be any filtration of the Borel sigma-field on $E$. Thus, the range of this section is not limited to stochastic analysis, and encompasses purely analytic frameworks, such as $E=S=\mathbb{R}$ with suitable filtrations (see Example~\ref{example1}). Lemma~\ref{consistentreg}, which is the key result of this section (resp. Proposition~\ref{characsmn}), characterizes {causal transference plans} (resp. {couplings}), in the particular case where $(\mathcal{B}_t(S))$ is regular. The assumptions on filtrations allow some communication of the topological properties of the underlying Polish space $S$, to the topology of several sets of transference plans  (Theorem~\ref{theoremcompact2}), through Lemma~\ref{consistentreg} ; the latter encapsulates  the effects of the filtration's regularity.  This yields the \textit{primal attainment}, and the existence of \textit{dual formulations}, to several problems of \textit{causal optimization} (Corollary~\ref{ToSEP}). Dual formulations are derived from Theorem~\ref{thmduality}, by substituting for $\mathcal{H}_{\eta,\nu}$, the set $\mathcal{H}_\eta$, which does not depend on the second marginal.

\subsection{Assumptions on $(\mathcal{B}_t(S))$} In view of applications, the following assumption on filtrations on $S$ is crucial. The notation of Definition~\ref{definitionconsist} is inherited from \cite{I-W} p.165, which emphasizes these properties. \begin{definition}
\label{definitionconsist}
We say that a filtration $(\mathcal{B}_t(S))_{t\in I}$, of the Borel sigma field of a Polish space $S$, is a regular filtration, if there exists a family $(\rho_t)_{t\in I}$  of maps which meets the following assumptions \begin{enumerate}[(i)]
\item  For all $t\in I$ the map $\rho_t : S\to S$ is continuous, and $$\mathcal{B}_t(S) = \rho_t^{-1}(\mathcal{B}(S)).$$
\item  $\rho_s\circ \rho_t = \rho_{s\wedge t},$ for all $s,t\in I$, where $s\wedge t := min(s,t)$.
\item  There exists some $t_f\in I$, such that $\rho_{t_f}= I_S ;$ the identity map on $S$.
\end{enumerate}
In this case we say that $(\mathcal{B}_t(S))$ is induced by the consistent family of continuous map $(\rho_t)_{t\in I}$. Moreover, in this case, if $I:=[0,1]$, and if the map $\rho: (t,y) \in [0,1]\times S \to \rho_t(y)\in S$ is continuous, then we say that $(\mathcal{B}_t(S))$ is completely regular.
\end{definition}
Assume that $S$ is endowed with a \textit{regular filtration} $(\mathcal{B}_t(S))_{t\in I}$, and that $(\Omega, \mathcal{A},  \mathcal{P})$ is a complete probability space, with a $\mathcal{P}-$complete filtration $(\mathcal{A}_t)_{t\in I}$, not necessarily {right-continuous}. To shorten subsequent statements, we say that a $\mathcal{A}/\mathcal{B}(S)-$measurable map $\widetilde{X} : \Omega \to S$, is an $(\mathcal{A}_t)-$adapted map, if $$ \widetilde{X}^{-1}(\mathcal{B}_t(S)) \subset \mathcal{A}_t,$$ for all $t\in I ;$ such maps are usually handled in stochastic analysis, to perform transformations of measures.  If $I=[0,T]$, for some $T>0$, then by following \cite{JACOD}, for convenience of notation,  by a \textit{complete stochastic basis} $(\Omega, \mathcal{A}, (\mathcal{A}_t)_{t\in I}, \mathcal{P})$, we mean a complete probability space $(\Omega, \mathcal{A}, \mathcal{P}),$ together with a filtration $(\mathcal{A}_t)$, such that $\mathcal{A}_t\subset \mathcal{A},$ for all $t\in I$, which further satisfies the usual conditions (i.e. it is $\mathcal{P}-$complete, and \textit{right-continuous} ;  $\mathcal{A}_t= \mathcal{A}_{t+}$, for all $t\in I$).
\begin{example}
\label{example1}
\begin{enumerate}[(i)]
\item Setting $I:=\mathbb{R}\cup\{+\infty\}$, we have  $\mathcal{B}_t(\mathbb{R}) = {\rho^\mathbb{R}_t}^{-1}(\mathcal{B}(\mathbb{R})),$ where  \begin{equation} \label{btrdef} \mathcal{B}_t(\mathbb{R}) := \sigma((-\infty, a], a <t),\end{equation} and where $\rho_t^\mathbb{R}$ denotes the continuous map  $$\rho_t^\mathbb{R} : s\in \mathbb{R} \to s\wedge t \in \mathbb{R},$$ for all $t\in I$. Thus, $(\mathcal{B}_t(\mathbb{R}))_{t\in I}$ is a regular filtration.
\item Given $T>0$,  $\mathcal{B}([0,T])$ denotes the usual sigma-field on $[0,T]$, which is the trace of $\mathcal{B}(\mathbb{R})$ on $[0,T]$.  We have, $\mathcal{B}_t([0,T]) = {\rho_t^{[0,T]}}^{-1}\left(\mathcal{B}([0,T])\right),$ where $\mathcal{B}_t([0,T])$ denotes the trace of $\mathcal{B}_t(\mathbb{R})$ on $[0,T]$, and where $$\rho_t^{[0,T]} : s\in [0,T] \to s\wedge t \in [0,T],$$ for all $t\in[0,T]$. Take $(\Omega, \mathcal{A}, (\mathcal{A}_t)_{t\in[0,T]}, \mathcal{P})$, a complete probability space, with a $\mathcal{P}-$complete filtration $(\mathcal{A}_t)$, such that $\mathcal{A}_t\subset \mathcal{A},$ for all $t\in[0,T]$.  A $\mathcal{A}/ \mathcal{B}([0,T])$-measurable map $\tau : \omega \in \Omega \to \tau(\omega) \in[0,T]$ is $(\mathcal{A}_t)-$adapted, if and only if,   $\{\tau < t\} \in \mathcal{A}_t,$ for all $t\in [0,T]$ iff $\tau$ is an $(\mathcal{A}_{t+})-$\textit{stopping time} (bounded by $T$). In particular, if $(\Omega, \mathcal{A}, (\mathcal{A}_t)_{t\in[0,T]}, \mathcal{P})$ is a complete stochastic basis, then $\tau$ is an $(\mathcal{A}_t)-$stopping time, if and only if, $\tau$ is an $(\mathcal{A}_t)-$adapted map. By considering $[0,\infty]$ with its usual compact Polish structure, this extends to unbounded stopping times.
\item The space $\mathbb{R}^n$, for some $n\in \mathbb{N}$, with $I:= \{1,...,n\}$ is naturally endowed with the filtration $(\mathcal{B}_t(\mathbb{R}^n))_{t\in I}$, which is defined by $
\mathcal{B}_t(\mathbb{R}^n) := \sigma(\pi_i, i\leq t),$ for $t\in I$, where $\pi_i : (x_1,...,x_n)\in \mathbb{R}^n \to x_i \in \mathbb{R}$ denotes the canonical projection on the $i-th$ 
component. Defining, the continuous map  $$\rho^n_t : x\in \mathbb{R}^n \to (x_1,...,x_{t},0,...,0) \in \mathbb{R}^n,$$ we obtain $\mathcal{B}_t(\mathbb{R}^n)= 
({\rho_t^n})^{-1}(\mathcal{B}(\mathbb{R}^n)),$ for all $t\in I$. If $X :\Omega \to \mathbb{R}^n$ denotes a measurable map, on a complete probability space $(\Omega, \mathcal{A}, \mathcal{P})$, with a $\mathcal{P}-$complete filtration $(\mathcal{A}_t)_{t\in I}$ of $\mathcal{A}$, define $X_t:= \pi_t\circ X,$ for all $t\in I$. Thus,  $(X_t)_{t\in I}$ is a stochastic 
process, on this probability space, and  $X^{-1}(\mathcal{B}_t(\mathbb{R}^n)) \subset \mathcal{A}_t,$ for all $t\in I$, if and only if, $X_t$ is $\mathcal{A}_t-$measurable, for all $t\in I$.
\item On $[0,T]^n$ (cartesian product), for some $n\in \mathbb{N}$, the filtration $(\mathcal{B}_t([0,T]^n))_{t\in [0,T]}$ defined by $\mathcal{B}_t([0,T]^n) := \mathcal{B}_t([0,T])\otimes...\otimes \mathcal{B}_t([0,T]),$ for all or all $t\in[0,T]$, is a \textit{regular filtration}. Indeed,  $$\widetilde{\rho}_t : (s_1,...,s_n) \in [0,T]^n \to (s_1 \wedge t, ... s_n\wedge t) \in [0,T]^n$$ yields $(\mathcal{B}_t([0,T]^n)) = (\widetilde{\rho}_t^{-1}(\mathcal{B}([0,T]^n))).$ Define  $$U:= \tau_1 \times... \times \tau_n : \omega \in \Omega \to (\tau_1(\omega),..., \tau_n(\omega)) \in [0,1]^n,$$ where $\tau_i : \Omega \to [0,T]$ is some $\mathcal{A}/ \mathcal{B}([0,T])-$measurable map on a \textit{complete stochastic basis} $(\Omega, \mathcal{A}$ $, (\mathcal{A}_t)_{t\in[0,T]}, \mathcal{P})$. All the $(\tau_i)$ are $(\mathcal{A}_t)$-stopping times (bounded by $T$), if and only if, $U$ is a $(\mathcal{A}_t)-$adapted map.  It extends  to unbounded stopping times using the remark in $(ii)$.

\item  In this paper,  $W$ denotes the space $C([0,1], \mathbb{R}^d)$, of  $\mathbb{R}^d-$valued continuous maps on $[0,1]$, endowed with the norm of uniform convergence. In stochastic calculus, each $\omega\in W$  models a sample path, of a $\mathbb{R}^d-$valued continuous process, on the interval of time $[0,1]$. The natural filtration $(\mathcal{B}_t^0(W))_{t\in[0,1]}$ is defined by setting $\mathcal{B}_t^0(W) := \sigma(W_s, s\leq t),$ for $t\in [0,1]$, where $W_t$ is the map $$W_t : \omega\in W \to W_t(\omega):= \omega(t) \in \mathbb{R}^d,$$ for  $t\in[0,1]$.  Given $\eta\in \mathcal{P}_W$, $(W_t)_{t\in [0,1]}$ defines the so-called evaluation process on the probability space $(W,\mathcal{B}(W)^\eta,\eta)$.  Setting $$\rho_t^W : \omega \in W \to \omega_{.\wedge t} \in W,$$ yields $\mathcal{B}_t^0(W) = {\rho_t^W}^{-1}(\mathcal{B}(W)),$ for all $t\in[0,1]$, so that  $(\mathcal{B}_t^0(W))_{t\in[0,1]}
$ is a regular filtration. Moreover,  it is completely regular. Further considering a complete stochastic basis $(\Omega, \mathcal{A},  (\mathcal{A}_t)_{t\in [0,1]}, \mathcal{P})$, and a measurable map $X : \Omega \to W$, define $$X_t : \omega \in 
\Omega \to (W_t\circ X)(\omega):= W_t(X(\omega))\in \mathbb{R}^d,$$ for all $t\in[0,1]$. We call $(X_t)_{t\in[0,1]}$ the process associated to $X$ on $(\Omega,\mathcal
{A}, \mathcal{P})$. The process $(X_s)_{s\in [0,1]}$ is $(\mathcal{A}_t)$-adapted, in the usual acceptation of stochastic 
calculus (i.e. $X_t$ is $\mathcal{A}_t$-measurable, for all $t\in [0,1]$), if and only if, $X$ is an $(\mathcal{A}_t)-$adapted map. 

\item Consider the \textit{complete stochastic basis}  of $(v)$, with the same map $X$, and the associated process $(X_s)_{s\in [0,1]}$.  We still define $S=W$, but we set $(\mathcal{B}_t(W)):= (\mathcal{B}_{t+}^0(W))$. It is not a \textit{regular filtration}, however since $(\mathcal{A}_t)$ satisfies the usual conditions (in particular it is \textit{right-continuous}), it defines the same $(\mathcal{A}_t)-$adapted maps as the filtration of $(v)$.  Moreover, setting  $ \mathcal{G}_t^X:= X^{-1}(\mathcal{B}_{t}(W))^{\mathcal{P}}$, $(\mathcal{G}_t^X)$ is the usual augmentation, of the filtration generated by the stochastic process  $(X_s)_{s\in [0,1]}$,  in the usual acceptation of stochastic calculus. \end{enumerate}
\end{example}

\subsection{Causal transference plans, for regular filtrations on $S$}  In Lemma~\ref{consistentreg} below, the equivalence of $(i)$ and $(ii)$ states that, when the filtration on $S$ is regular, causal transference plans are exactly those transference plans, which satisfy the constraint on joint laws introduced in the proof of the Yamada-Watanabe criterion (see \cite{I-W} p.164-166). The equivalence of $(ii)$ and $(iii)$ is trivial, and follows as a particular case of well known  results (see \cite{Jacod2}, for instance), we state it precisely, and provide a concise proof, for the sake of completeness. \begin{lemma}
\label{consistentreg}
 Given $\eta \in \mathcal{P}_E$, $\nu \in \mathcal{P}_S$ and $\gamma \in \Pi(\eta,\nu)$, by further assuming that $(\mathcal{B}_t(S))$ is a {regular filtration}, for all $t\in I$, $\mathcal{G}_t^\gamma$ (see Definition~\ref{defmain}) is also the $\eta-$completion, of the smallest sigma-field on $E$, such that for all $A\in \mathcal{B}_t(S)$, the map $$\phi_A: x\in E\to \phi_A(x):= \Theta^x(A) \in \mathbb{R}$$ is measurable, $\Theta$ denoting any kernel associated to $\gamma$ by~(\ref{kernot}). In particular, the following assertions  are equivalent \begin{enumerate}[(i)] \item $\gamma \in \Pi_c(\eta,\nu)$
 \item $\gamma$ satisfies the Yamada-Watanabe constraint. That is, for all $t\in I$, and $A\in \mathcal{B}_t(S)$,  the map $\phi_A$ is $\mathcal{B}_t(E)^\eta-$measurable.
 \item For all $f\in \mathcal{L}^\infty(\eta)$ $$E_\gamma[f\circ\pi | (\mathcal{B}_t(E)\otimes \mathcal{B}_t(S))^\gamma] =  E_\eta[f| \mathcal{B}_t(E)^\eta] \circ \pi \ \gamma-a.s., \ for \ all \ t\in I.$$
 \end{enumerate}
\end{lemma}
\nproof Given $t\in I$, denote by $\mathcal{H}_t$, the $\eta-$completion of the smallest sigma-field on $E$,  such that, the map $\phi_A$  is measurable, for all $A\in \mathcal{B}_t(S)$. By definition $\mathcal{G}^\gamma_t \subset \mathcal{H}_t$ ; we want yo obtain the converse inclusion. By definition, it is equivalent to prove, that the map $\phi_A$ is $\mathcal{G}_t^\gamma-$measurable, for all $A\in \mathcal{B}_t(S)$. For any element $L$ of the set $L^\infty_{1,+}:= \left\{ L \in L^\infty(\eta) \middle| \eta-a.s. L>0, E_\eta[L] = 1 \right\},$ by Fubini's theorem (see \cite{DM} p.31-33),  we define two probabilities on the measurable space $(S,\mathcal{B}_t(S))$, by setting $\widetilde{\mathcal{P}}_{t,L} := \int_E d\eta(x) L(x) \Theta^x,$  and  $\widetilde{\mathcal{Q}}_{t,L} := \int_E d\eta(x) E_\eta[L | \mathcal{G}^\gamma_t] \Theta^x.$  Denote by $(\rho_s)$, a family of continuous maps, which induces $(\mathcal{B}_s(S))$, in the acceptation of Definition~\ref{definitionconsist}. Since $\rho_t$ is $\mathcal{B}_t(S) / \mathcal{B}(S)-$measurable, setting $\mathcal{P}_{t,L} := {\rho_t}_\star \widetilde{\mathcal{P}}_{t,L},$ and $\mathcal{Q}_{t,L} := {\rho_t}_\star \widetilde{\mathcal{Q}}_{t,L}$, we obtain two Borel probabilities on $(S,\mathcal{B}(S))$, for all $L\in L^\infty_{1,+}$. By linearity,  $\mathcal{H}_t  \subset \mathcal{G}^\gamma_t$, is equivalent to
$\mathcal{P}_{t,L}= \mathcal{Q}_{t,L}$, for all $L\in  L^\infty_{1,+}$. On the other hand, from the definitions, the map $\phi_A$ is $\mathcal{G}^\gamma_t-$measurable, for all $A\in \mathcal{B}_t(S)$ of $\nu-$continuity. Whence, we obtain \begin{equation}  \label{equalm} \mathcal{P}_{t,L}(B) = \mathcal{Q}_{t,L}(B),   \end{equation} for all $L\in  L^\infty_{1,+}$, and $B\in \mathcal{B}(S)$, such that $\nu(\partial \rho_t^{-1}(B))=0$. Since $\rho_t$ is continuous, and since, for all $L\in  L^\infty_{1,+}$, ${\rho_t}_\star \nu\sim\mathcal{P}_{t,L}$ (i.e. {equivalent}),~(\ref{equalm}) holds, for all $L\in  L^\infty_{1,+}$, and $B\in \mathcal{B}(S)$, such that $\mathcal{P}_{t,L}(\partial B) = 0$. Thus, since $\mathcal{P}_{t,L} \in \mathcal{P}_E$, the Portmanteau theorem (see \cite{SV}) yields $\mathcal{P}_{t,L}= \mathcal{Q}_{t,L}$, for all $L \in  L^\infty_{1,+}$ ; this proves that $\mathcal{G}_t^\gamma = \mathcal{H}_t$.  In particular $(i)$ is equivalent to $(ii)$. On the other hand, for all $t\in I$,  and $f\in \mathcal{L}^\infty(\eta)$, we have  $$E_\gamma\left[(f\circ \pi - E_\eta[f| \mathcal{B}_t(E)^\eta]\circ \pi) 1_A\circ \pi 1_B \circ \widetilde{\pi}\right] =E_\eta\left[ f1_A\left(\Theta(B) -E_\eta[\Theta(B)| \mathcal{B}_t(E)^\eta]\right)\right],$$ for all $A\in \mathcal{B}_t(E)$, $B\in \mathcal{B}_t(S)$. Therefore, the equivalence of $(ii)$ with $(iii)$ follows, by a monotone class argument. \nqed    

From Definition~\ref{defmain}, Lemma~\ref{consistentreg} (resp.  Lemma~\ref{consistentreg}, together with a monotone class argument),  easily entails the following Proposition~\ref{consistentreg2} (resp. Proposition~\ref{characsmn}), which characterizes deterministic causal transference plans (resp. causal couplings), under this further assumptions on $(\mathcal{B}_t(S))$. 
  \begin{proposition}
\label{consistentreg2}
 Given $\eta \in \mathcal{P}_E$, $\nu \in \mathcal{P}_S$, $U \in \mathcal{R}(\eta,\nu)$, and further assuming that $(\mathcal{B}_t(S))$ is a {regular filtration}, 
for $\mathcal{G} \subset \mathcal{B}(S)$ a sigma field, we define $U^{-1}(\mathcal{G}):= f^{-1}(\mathcal{G})^\eta,$ for any $\mathcal{B}(E)^\eta/\mathcal{B}(S)-$ 
measurable map $f: E\to S$, whose associated $\eta-$equivalence class of maps is $U$. Then, we have $$\mathcal
{G}_t^{j(U)}= U^{-1}(\mathcal{B}_t(S)),$$ for all $t\in I$,  $j$ denoting the map given by~(\ref{jdef}).  In particular, $U\in \mathcal{R}_a(\eta, \nu)$ if and only if $U^{-1}(\mathcal{B}_t(S))\subset \mathcal{B}_t(E)^\eta,$ for all $t\in I$. \nqed
\end{proposition}
  \begin{proposition}
  \label{characsmn}
 On a complete probability space $(\Omega,\mathcal{A},\mathcal{P})$, let $Y : \Omega \to E$,  and  $Z:\Omega \to S$, be two $\mathcal{A} / \mathcal{B}(E)-$ (resp. $\mathcal{A} / \mathcal{B}(S)-$) measurable maps. Further define  $(\mathcal{G}^Y_t):=(Y^{-1}(\mathcal{B}_t(E)))^{\mathcal{P}}$ (resp. $(\mathcal{G}_t^Z):= (Z^{-1}(\mathcal{B}_t(S)))^{\mathcal{P}} $ ), the completed filtration generated by $Z$ (resp. by $Y$). Further assuming that $(\mathcal{B}_t(S))$ is a {regular filtration}, the following assertions are equivalent 
  \begin{enumerate}[(i)]
  \item $(Y\times Z)_\star \mathcal{P} \in \Pi_c( Y_\star \mathcal{P}, Z_\star \mathcal{P})$,  i.e. $(Y,Z)$ is a {causal coupling}.
  \item We have  $$\mathcal{P}\left( Z\in C \middle| \mathcal{G}^Y_t\right) =\mathcal{P}\left( Z\in C \middle| \sigma(Y)^{\mathcal{P}}\right) \ \mathcal{P}-a.s.,$$ for any $C\in \mathcal{B}_t(S)$, and $t\in I $, $\sigma(Y)^{\mathcal{P}}$ denoting the $\mathcal{P}-$completion of $Y^{-1}(\mathcal{B}(E))$.
  \item  We have $$E_{\mathcal{P}}[f\circ Y | \mathcal{G}_t^Y] =  E_{\mathcal{P}}[f\circ Y | \sigma( \mathcal{G}_t^Y \cup \mathcal{G}_t^Z) ] \ \mathcal{P}-a.s.,$$ for any $t\in I $, and any $f\in L^1(Y_\star \mathcal{P})$.
  \end{enumerate} \nqed
 \end{proposition}
\begin{remarkk}(Stability by pullback)\label{pullbackcaustp} Consider $E$, $S$, $Z$ three Polish spaces, endowed with filtrations of their Borel-sigma fields, those of $S$ and $Z$ being  \textit{regular filtrations}. Take $\eta\in \mathcal{P}_E$, $\nu\in \mathcal{P}_S$, $\widetilde{\nu} \in \mathcal{P}_Z$. For $\gamma \in \Pi(\eta,\nu)$ (resp. $\widetilde{\gamma} \in \Pi(\nu,\widetilde{\nu}))$, we denote by $\Theta_\gamma$ (resp. by $
\Theta_{\widetilde{\gamma}}$) an associated \textit{kernel} to $\gamma$ (resp. to $\widetilde{\gamma}$) by~(\ref{kernot}). Setting $\Theta_{\widetilde{\gamma}\circ \gamma}: (x,A)\in E \times \mathcal{B}(Z) \to \Theta_{\widetilde{\gamma}\circ \gamma}^x(B) = \int_S \Theta_\gamma^x(dy) \Theta_{\widetilde{\gamma}}^y(B) \in \mathbb{R},$ by Fubini's theorem, it defines uniquely a probability, which we denote by $\widetilde{\gamma}\circ \gamma \in \Pi(\eta,\widetilde{\nu}),$  such that $\widetilde{\gamma}\circ \gamma(A\times B) = E_{\eta}[1_A \Theta_{\widetilde{\gamma}\circ \gamma}(B) ],$ for $A\in \mathcal{B}(E)$, $B\in \mathcal{B}(Z)$.  Under our assumptions on marginals, it does not depend on the version of the \textit{kernels}. We have $j(V)\circ j(U)= j(V\circ U),$ for all $U\in \mathcal{R}(\eta,\nu)$, $V\in \mathcal{R}(\nu,\widetilde{\nu})$, where $j$ is the map defined by~(\ref{jdef}), and where in the right hand term , $\circ$ denotes the pullbacks of morphisms of probability spaces (see \cite{MAL2} p.156). Further assuming that $\gamma \in \Pi_c(\eta,\nu),$ and $\widetilde{\gamma} \in \Pi_c(\nu,\widetilde{\nu})$, Fubini's theorem (see \cite{DM}) yields,  $\widetilde{\gamma}\circ\gamma \in \Pi_c(\eta,\widetilde{\nu}).$ This allow to define distance on the Wasserstein space, based on the symmetric counterpart of causal couplings defined in Section~\ref{4}, with similar proof as in the unconstrained case. \end{remarkk}  
\subsection{Further topological properties, and the associated optimization problems}
\begin{theorem}
 \label{theoremcompact2}
 For $\eta \in \mathcal{P}_E$, define $\mathcal{P}^S_c(\eta):= \cup_{\nu\in \mathcal{P}_S} \Pi_c(\eta,\nu),$ and $L^0_a(\eta,S):= \cup_{\nu\in \mathcal{P}_S}  \mathcal{R}_a(\eta,\nu) ;$ see  Definition~\ref{defmain}. Further assuming that $S$ is endowed with a {regular filtration}, these sets have the following properties 
\begin{enumerate}[(i)]
\item $\mathcal{P}^S_c(\eta)$ is convex and closed in $\mathcal{P}_{E\times S}$, for the topology of weak convergence in measure. In particular, if $S$ is further assumed to be compact for its Polish topology, then $\mathcal{P}^S_c(\eta)$ is compact.
\item  $L^0_a(\eta, S)$ is a closed subset of $L^0(\eta, S)$ (see Section~\ref{2}), for the topology of convergence in probability. 
\item Let $(\gamma_n)$ be a sequence of elements of $\mathcal{P}_c^S(\eta)$. Further assuming that $(\widetilde{\pi}_\star \gamma_n)$ is tight, there exists a $\gamma \in \mathcal{P}^S_c(\eta)$, and a subsequence $(\gamma_{k(n)})$ of $(\gamma_n)$, which converges weakly in measure to $\gamma$.
\end{enumerate}  \end{theorem}  \nproof The convexity in $(i)$  is trivial.  By continuity of the map $j$ defined by~(\ref{jdef}), $(ii)$ follows from $(i)$.  Whenever $(\widetilde{\pi}_\star \gamma_n)$ is tight, $\cup_{n} \Pi(\eta, \widetilde{\pi}_\star \gamma_n)$ is tight (see \cite{Vil2} p.45). Thus, we extract a subsequence  $(\gamma_{k(n)}) \subset \mathcal{P}^S_c(\eta)$, which converges to a $\gamma$. By continuity of $\pi_\star$, $\gamma \in \mathcal{P}^S(\eta)$.  Hence,  $(iii)$ follows from $(i)$. On the other hand, together with Lemma~\ref{consistentreg}, it follows, similarly to the proof of Theorem~\ref{theoremcompact}, that $\mathcal{P}^S_c(\eta)$ is closed. If $S$ is compact, by Prohorov's criterion, the inner regularity of $\eta$ implies the compactness of $\mathcal{P}_S(\eta)$. Thus, the closed subset $\mathcal{P}_c^S(\eta)$ of $\mathcal{P}^S(\eta)$ is compact. \nqed

We obtain the following relaxation to Theorem~\ref{thmduality} ; recall  that   
 $\mathcal{H}_\eta$ is introduced in Definition~\ref{heta}
\begin{corollary}
\label{thmduality1}
Given a non negative l.s.c. map $c: E\times S \to \mathbb{R}\cup \{+\infty\}$, we further assume that $(\mathcal{B}_t(S))$ is a regular filtration. For any $\eta\in \mathcal{P}_E,\nu \in \mathcal{P}_S$, we have $D^{R}_{\eta,\nu} = P_{\eta,\nu},$ where $$P_{\eta,\nu}:= \inf\left(\left\{ \int_{E\times S} c(x,y)d\gamma(x,y) \middle| \gamma \in \Pi_c(\eta,\nu) \right\} \right),$$ and where $$D^R_{\eta,\nu} =\sup\left(\left\{\int_S g(y) d\nu(y) \middle| (g,h) \in C_b(S)\times\mathcal{H}_\eta : g\circ\widetilde{\pi} + h \leq c\right\}\right).$$ 
\end{corollary}
\nproof
Since $\mathcal{H}_{\eta,\nu} \subset \mathcal{H}_\eta$, $P_{\eta,\nu}  \leq D^R_{\eta,\nu}$ follows from Theorem~\ref{thmduality}. On the other hand, similarly to the proof of $(ii)$ of Lemma~\ref{pmropo3}, Lemma~\ref{consistentreg} easily implies  $E_\gamma[h]=0,$ for all $h\in \mathcal{H}_\eta$, and $\gamma \in \Pi_c(\eta,\nu)$. Whence, the converse inequality follows, from the definition.
\nqed

Recall that, given $\eta\in \mathcal{P}_E$,  $\mathcal{P}_c^S(\eta)$ is defined by Theorem~\ref{theoremcompact2}. 
\begin{corollary}
\label{ToSEP}
 Given a non-negative l.s.c. map $c : E\times S  \to \mathbb{R} \cup\{+ \infty\},$ let   $$e : E\times S \to Z,$$ be a continuous map with values in a third Polish space $Z$. For $\eta \in \mathcal{P}_E,$ and $\widetilde{\mu} \in \mathcal{P}_Z$, define the following {primal problems}  \begin{equation} \label{absos} P_\eta := \inf\left(\left\{ \int_{E\times S} c(x,y) d\gamma(x,y) \middle| \gamma \in \mathcal{P}_c^S(\eta)\right\}\right), \end{equation}  and  \begin{equation} \label{sepabstract} P^e_{\eta,\widetilde{\mu}} := \inf\left(\left\{ \int_{E\times S} c(x,y) d\gamma(x,y)\middle| \gamma \in \mathcal{P}_c^S(\eta) , e_\star \gamma =\widetilde{\mu}\right\}\right),\end{equation} and the related dual problems $$D_\eta := \sup\left(\left\{ a \middle|  (a, h) \in \mathbb{R}\times \mathcal{H}_\eta \ :  \ a + h(x,y)  \leq c(x,y),  \forall (x,y)\in E\times S \right\}\right),$$ and $$D^e_{\eta,\widetilde{\mu}} := \sup\left(\left\{\int_Z  g(z) d\widetilde{\mu}(z)  \middle| (g,h) \in C_b(Z)\times \mathcal{H}_\eta \ : \  h+  g\circ e \leq c \right\}\right),$$  where the set $\mathcal{P}_c^S(\eta)$ was defined in Theorem~\ref{theoremcompact2}. Further assuming that $S$ is compact for its Polish topology, and that $S$ is endowed with a regular filtration, the infimum in~(\ref{absos}) is always  attained (resp. if $P^e_{\eta,\widetilde{\mu}} <\infty$ the infimum in~(\ref{sepabstract}) is attained), and $D_\eta= P_\eta$ (resp. by assuming that $P^e_{\eta,\widetilde{\mu}}<\infty$, $D^e_{\eta,\widetilde{\mu}}= P^e_{\eta,\widetilde{\mu}}$).
\end{corollary}
\nproof Since $S$ is compact, by Theorem~\ref{theoremcompact2}, $\mathcal{P}_c^S(\eta)$ is compact. On the other hand, by the Portmanteau theorem (see \cite{STROOCK}), the assumptions on $c$ entail the lower semi-continuity of $\gamma \in \mathcal{P}_{E\times S} \to E_\gamma[c]\in \mathbb{R}\cup\{+\infty\}$.  Thus, the infimum of $P_\eta$ is always attained. Define $$\phi : \gamma \in \mathcal{P}_{E\times S} \to e_\star \gamma \in 
\mathcal{P}_Z,$$ so that $$P^e_{\eta,\widetilde{\mu}} = \sup\left(\left\{E_\gamma[c] \middle| \gamma \in \mathcal{P}_c^S(\eta)\cap \phi^{-1}(\{\widetilde{\mu} \}) \right\}\right) .$$Since $e$ is continuous, $\phi$ is continuous for the topology of weak convergence in measure. Thus,  $\mathcal{P}_c^S(\eta)\cap \phi^{-1}(\{\widetilde{\mu}\})$ is compact, as a 
closed subset, of the compact set $\mathcal{P}_c^S(\eta)$. Similarly, the hypothesis on $c$ entail that if $P^e_{\eta,\widetilde{\mu}} <\infty$, then the infimum is attained ; however  $\mathcal{P}_c^S(\eta)\cap \phi^{-1}(\{\widetilde{\mu}\})$ can be empty if $P^e_{\eta,\widetilde{\mu}} =\infty$.  By Lemma~\ref{consistentreg}, similarly to the proof of $(ii)$ of Lemma~\ref{pmropo3}, we obtain $E_\gamma[h]=0$, for all $h\in \mathcal{H}_\eta$, $\gamma \in \mathcal{P}_c^S(\eta)$. Thus, similarly to the proof of Corollary~\ref{thmduality1}, we obtain $D_\eta\leq P_\eta$ (resp. $D^e_{\eta,\widetilde{\mu}}\leq P^e_{\eta,\widetilde{\mu}}$). Using Moreau-Yosida's approximation, it is enough to prove the converse inequalities under the assumption that $c\in C_b(E\times S)$. Under this assumption, Corollary~\ref{thmduality1} yields $$P_{\eta}= \inf_{\nu\in \mathcal{P}_S} \inf_{\gamma\in \Pi_c(\eta, \nu)} E_\gamma[c] = \inf_{\nu\in \mathcal{P}_S} \sup_{k \in \mathcal{K}} F(\nu, k),$$ where, $\mathcal{K}$ denotes the convex set, of $k\in C_b(S)$, such that there exists a $h\in \mathcal{H}_\eta$, which satisfies  $h + k\circ \widetilde{\pi} \leq c,$ and where  $F : (\nu, k)\in \mathcal{P}_S\times \mathcal{K}  \to F(\nu,k):= E_\nu[k] \in \mathbb{R}.$ Since $S$ is compact,  $\mathcal{P}_S$ is compact ; on the other hand, for all $k\in \mathcal{K}\subset \mathcal{C}_b(S)$ (resp. $\nu \in \mathcal{P}_S$) , $F(.,k)$ is a convex continuous map (resp. a concave map). Thus, by a classical min-max theorem (for instance Theorem 2 of \cite{FAN}), \begin{equation} \label{minmaxlfmr} P_\eta = \sup_{k\in \mathcal{K}}  \inf_{\nu\in \mathcal{P}_S} E_\nu[k].\end{equation} For $k\in \mathcal{K}$, define $a_k:=  \inf_{\nu\in \mathcal{P}_S} E_\nu[k]$.  Taking $\nu= \delta^{Dirac}_y$, we obtain $h(x,y) +a_k \leq h(x,y) + k(y) \leq c(x,y),$  for all $(x,y)\in E\times S$. Therefore, $P_\eta = \sup_{k\in \mathcal{K}} a_k \leq  D_\eta,$ follows from~(\ref{minmaxlfmr}). Moreover we have $$ P^e_{\eta,\widetilde{\mu}} = \inf_{\gamma \in \mathcal{P}_c^S(\eta)} \sup_{g\in C_b(Z)} G(\gamma, g),$$ where  $G(\gamma, g) := E_\gamma[c - g\circ e + E_{\widetilde{\mu}}[g]]$. Since $c\in C_b(E\times S)$, and $g\in C_b(Z)$,  we obtain similarly $$ P^e_{\eta,\widetilde{\mu}} = \sup_{g\in C_b(Z)} (E_{\widetilde{\mu}}[g] +  \inf_{\gamma \in \mathcal{P}_c^S(\eta)}E_\gamma[c - g\circ e ]),$$  from the classical min-max theorem of \cite{FAN}.   Finally, substituting the continuous and bounded cost map $c-g\circ e$ for $c$, in the definition of the primal (resp. dual) problem $P_\eta$ (resp. $D_\eta$), we already proved that $D_\eta= P_\eta$ holds. Whence, we conclude that $$ P^e_{\eta,\widetilde{\mu}} = \sup_{g\in C_b(Z)} \sup_{a \in \mathbb{R}, \exists h\in \mathcal{H}_\eta : a + h \leq c-g\circ e}  (E_{\widetilde{\mu}}[g]  + a) \leq D^e_{\eta,\widetilde{\mu}}.$$ \nqed
\begin{remarkk} Under the assumptions of Corollary~\ref{thmduality1} and of Corollary~\ref{ToSEP}, by Lemma~\ref{consistentreg}, a transference plan is causal iff it induces an extension of filtered probability space (see the Introduction). Whence, taking suitable conditional expectations of elements of $\mathcal{H}_\eta$, which express the causal constraints through the penalization of the cost, the latter Corollaries can be easily formulated, in the Polish framework, to involve martingales. Moreover, in the particular case  $S=[0,1]$, with the filtration of Example~\ref{example1}, $\eta$ denoting laws of solutions to suitable stochastic differential equations (see Section~\ref{4}, and Section~\ref{5} below), by the straightforward representation provided in Section~\ref{4} and Section~\ref{5}, below, the dual formula of~\cite{BHC} can be obtained, and extended to non-Markovian frameworks. Finally, the same approach yields similar dual formulas, involving martingales, for optimal stopping problems ; to avoid an overlap with \cite{BHC}, those are not investigated in this paper.   \end{remarkk}

\section{Causal coupling plans in stochastic frameworks, and their symmetric counterparts}
\label{4}
In Section~\ref{3}, $(\mathcal{B}_t(S))$ was assumed to be regular ; none assumption was required on $(\mathcal{B}_t(E))$. Subsequent applications to stochastic calculus arise by suitable choices of $E$, of its filtration, and of the first marginal $\eta \in \mathcal{P}_E$. The first subsection states, in Polish framework, the symmetric counterpart to causal coupling plans, which naturally appears in Section~\ref{6} below.  Then, for the sake of clarity, we focus on the case where $E$ is a space of $\mathbb{R}^d-$ valued continuous paths, and $I=[0,1]$. As it will be clear below, by $(ii)$ of Example~\ref{example1}, relating causal couplings to stopping times requires the right-continuity of  $(\mathcal{B}_t(E))^\eta$ ; the latter also yields the existence of c\`{a}d-l\`{a}g modifications to martingales, which is useful to avoid measurability issues, in applications to stochastic calculus. The latter may fail, in general, when $\eta$ is the law of a \textit{non-Markovian} continuous semi-martingale . This motivates the general framework, stated in the first subsection,  for all our subsequent applications. Under weak assumptions, Proposition~\ref{propja1} and Proposition~\ref{propja2},  \textit{characterize causal transference plans from  laws of solutions to stochastic differential equation}, to probabilities  on a Polish space, endowed with a regular filtration. The latter are essentially analytic reformulations of results of \cite{Jacod2}, and of the proof of the Yamada-Watanabe criterion (see \cite{I-W}). It expresses that, in this precise framework, \textit{causal transference plans are transference plans inducing extensions of filtered probability spaces} ; see the Introduction. Remark~\ref{meprl} is crucial : it indicates how to apply the two previous Propositions.

\subsection{The symmetric counterpart to causal couplings}\label{symsubsec}For the sake of simplicity, until the end of this subsection we assume that $E=S= Z$, for some Polish space $Z$, endowed with a regular filtration $(\mathcal{B}_t^0(Z))$, and we take $\eta,\nu\in \mathcal{P}_Z$. Under the above framework, we take $I:=[0,1]$,  $(\mathcal{B}_t(E)) =(\mathcal{B}^0_{t+}(Z)),$ and $(\mathcal{B}_t(S)) =(\mathcal{B}^0_{t}(Z))$. The continuous map  \begin{equation} \label{Rdef} R : (x,y)\in E\times S \to (y,x)\in E\times S,\end{equation} induces the continuous map  \begin{equation} \label{rstardef} R_\star : \gamma \in \mathcal{P}_{E\times S} \to R_\star \gamma \in \mathcal{P}_{E\times S}.\end{equation} Thus, define   \begin
{equation} \label{causcoupl} \Pi_{cs}(\eta,\nu):= \Pi_c(\eta,\nu) \cap R_\star^{-1}(\Pi_c(\nu,\eta)), \end{equation} which is the symmetric counterpart to $\Pi_c(\eta,\nu)$, and \begin{equation} \label{raz} \mathcal{R}_{as}(\eta, \nu) :=j^{-1}(\Pi_{cs}(\eta,\nu)),\end
{equation} the symmetric counterpart to $\mathcal{R}_a(\eta,\nu)$, $j$ denoting~(\ref{jdef}).  Since $j$ and $R_\star$ are continuous, from the definitions, Theorem~\ref
{theoremcompact} implies that $\Pi_{cs}(\eta,\nu)$ (resp. $\mathcal{R}_{as}(\eta,\nu))$ is compact (resp. closed) for the respective 
topologies.   Finally, given $\gamma \in \Pi(\eta, \nu)$, since $R_\star (R_\star \gamma) = (R\circ R)_\star  \gamma = \gamma,$ we obtain $\gamma\in \Pi_{cs}(\eta,\nu)$, if and only if,  $R_\star\gamma\in \Pi_{cs}(\nu,\eta)$ ; this will be used in Section~\ref{6}.

\subsection{Notation, and framework of subsequent applications}\label{soussection42} For the sake of clarity, in this paper the only path space we consider is $W:= C([0,1], \mathbb{R}^d)$,  the separable Banach space of $\mathbb{R}^d-$valued continuous 
paths on $[0,1]$, endowed with the norm $|.|_W$ of uniform convergence.  As in $(v)$ of~Example~\ref{example1}, the \textit{evaluation process} is defined by \begin{equation} \label{evaldefsde} W_t : \omega \in W \to W_t(\omega):= \omega(t)\in \mathbb{R}^d, \end{equation} and $(\mathcal{B}_t^0(W))_{t\in[0,1]}$  denotes its \textit{natural filtration}, which is defined by  \begin{equation} \label{evalprokf} \mathcal{B}_t^0(W)=\sigma(W_s, s\leq t),\end{equation} for $t\in[0,1]$.  Under the acceptation of Section~\ref{3}, it is a \textit{regular filtration}, with $(\mathcal{B}_t^0(W))= ({\rho_t^W}^{-1}(\mathcal{B}(W))),$ where \begin{equation} \label{rhotWdefcor} \rho_t^W: \omega \in W \to \omega_{.\wedge t}\in W,\end{equation} for all $t\in[0,1]$. In view of a unified framework for the applications, we further define $(\mathcal{B}_t(W))$ by  $$\mathcal{B}_t(W):= \mathcal{B}_{t+}^0(W):= \cap_{\epsilon>0} \mathcal{B}_{(t+\epsilon)\wedge 1}^0(W),$$ for $t\in[0,1]$. Given $\eta \in \mathcal{P}_W$, for 
consistency of notation with \cite{CRUZLAS}, $(\mathcal{F}_t^\eta)$ denotes the $\eta-$usual augmentation of $(\mathcal{B}_t^0(W))$ ; it is defined by \begin{equation} \label{emprol} 
\mathcal{F}_t^\eta = \mathcal{B}_{t}(W)^\eta, \end{equation} for all $t\in[0,1]$. To handle all  applications within the same framework, \textit{until the end of the paper, we take 
$I:=[0,1]$, $E=W,$ and we endow it with the filtration \begin{equation} \label{defilW}( \mathcal{B}_t(E))_{t\in[0,1]}:= (\mathcal{B}_{t}(W))_{t\in[0,1]}. \end{equation} }  Thus, $(\mathcal{B}_t(E))$ is not the canonical regular filtration $(\mathcal{B}_t^0(W))$ on $W$. Until the end of this section $S$ still denotes any Polish space endowed with a regular filtration $(\mathcal{B}_t(S))_{s\in[0,1]}$. In Section~\ref{5} (resp. Section~\ref{6}), we will take $S=[0,1]$ (resp. $S=W$) with the respective regular filtrations of Example~\ref{example1}. 
\begin{remarkk}(right-continuity of the filtration on $E$.) In the general framework of two different Polish spaces $E$ and $S$, denote by $(\mathcal{B}_t^0(E))$ (resp. by $(\mathcal{B}_t^0(S))$), a regular filtration on $E$ (resp. on $S$). Taking $(\mathcal{B}_t(E)):=(\mathcal{B}^0_{t+}(E))$  ensures the existence of \textit{c\`{a}d-l\`{a}g modifications} to any $(\mathcal{B}_t(E)^\eta)-$\textit{martingales} (see \cite{I-W}). On the other hand, to apply results of Section~\ref{3}, one will require a regular filtration on $S$ ; it  entails an apparent asymmetry. However, under these assumptions,  taking  $(\mathcal{B}_t(S)):=(\mathcal{B}_t^0(S))$, or $(\mathcal{B}_t(S)):=(\mathcal{B}_{t+}^0(S))$ define the same causal transference plans. 
  \end{remarkk}
\subsection{Causal transformations of law of SDE} Recall that $E=W$, $I=[0,1]$,  $(\mathcal{B}_t(E))$ is given by~(\ref{defilW}) ; Given $\eta\in \mathcal{P}_W$, $(\mathcal{F}_t^\eta)$ is defined by~(\ref{emprol}). Consider the stochastic differential equation \begin{equation} \label{SDEJ} dX_t = \sigma_t(X) dB_t + b_t(X)dt ; Law(X_0)=\eta_0, \end{equation} where $\eta_0\in \mathcal{P}_
 {\mathbb{R}^d}.$ We assume that the \textit{weak existence}, and that the \textit{weak uniqueness}, of solutions hold for~(\ref{SDEJ}). That is, we assume the existence of a \textit{weak solution} $(X,B)$ on some \textit{complete stochastic basis} $(\Omega, \mathcal{A}, (\mathcal{A}_t)_{t\in[0,1]}, \mathcal{P})$, and of a \textit{unique} $\eta\in \mathcal{P}_W$, such that  $X_\star \mathcal{P}=\eta$, for all $X$ which solves~(\ref{SDEJ}) for some Brownian, on a probability space $(\Omega,\mathcal{A}, \mathcal{P})$. For the sake of clarity, we \textit{further 
 assume} that \begin{equation} \label{intcond}E_{\eta}\left[\int_0^1 (\sigma_t \sigma_t^t)_{ij} ds \right]  <\infty, \end{equation}  for all $i,j=1,...,n$. The latter hypothesis can be dropped to obtain similar statement as below, up to standard localizations techniques, and by using {local martingales}.  From standard results on transformations of laws of continuous semi-martingales (see Proposition~3.1. of \cite{CRUZLAS}, for a statement with the same notation), on $(W, \mathcal{B}(W)^\eta, \eta)$, the evaluation process~(\ref{evaldefsde}) has the decomposition \begin{equation} \label{coordMRT} W_t = W_0 + 
 M_t^\eta + \int_0^t b_s(\omega)ds, \end{equation}  for some $(\mathcal{F}_t^\eta)-$martingales $(M_t^\eta)$, see \cite{I-W}. Moreover $$<{M^\eta}^i, {M^\eta}^j> =\int_0^. {{(\sigma_s.\sigma^t_s)}^{i,j}} ds \ \eta-a.s.,$$  $<{M^\eta}^i, {M^\eta}^j>$ denoting the \textit{predicable covariation 
 process} of $M^\eta$ (see \cite{JACOD}). In particular (see \cite{JACOD} p.174-175), we assume that $\eta$ has the \textit{predicable representation property} : for any c\`{a}d-l\`{a}g  $(\mathcal{F}_t^\eta)-${martingale} $(M_t)_{t\in[0,1]}$  on $(W,\mathcal{B}(W)^\eta, \eta)$, there exists a  $(\mathcal{F}_t^\eta)-${predicable} process $(A^M_t)$, such that we have $M := \int_0^. A^M_s dM_{s}^\eta, \ \eta-a.s.$. In particular, since $(\mathcal{F}_t^\eta)$, which is given by~(\ref{emprol}), satisfies the usual conditions, from Theorem 4, p.76 of \cite{DM2}, any $(\mathcal{F}_t^\eta)-$martingale on this complete probability space has a \textit{continuous modification}.  Under this assumption, we have the two following Proposition~\ref{propja1} and   Proposition~\ref{propja2}, which are essentially reformulations of results from \cite{Jacod2} and \cite{I-W}.  Recall that, until the end of the paper, we  take $E= W$, and that $(\mathcal{B}_t(E))$ is given by~(\ref{defilW}) ; in particular $(\mathcal{B}_t(E)^\eta) = (\mathcal{F}_t^\eta)$.  It concretizes the definition of {causal transference plans} in the statements below.
 
 \begin{proposition}
 \label{propja1}
Assume that the weak existence, and that the weak uniqueness of solutions hold for~(\ref{SDEJ}), and that $\eta$, the law of its solutions, satisfies~(\ref{intcond}). If $S$ is endowed with a regular filtration $(\mathcal{B}_t(S))_{t\in[0,1]}$, then given $(X,B)$, a weak solution to~(\ref{SDEJ}), on a {complete stochastic basis} $(\Omega, \mathcal{A}$, 
 $(\mathcal{A}_t)_{t\in[0,1]}, \mathcal{P})$, for all $\mathcal{A} / \mathcal{B}(S)-$measurable map $Y : \Omega \to S$ which satisfies  $$Y^{-1}(\mathcal{B}_t(S)) \subset \mathcal{A}_t, \ for \ all \ t\in [0,1],$$ the joint law  $\gamma := (X\times Y)_\star \mathcal{P}$ has the following properties :
  \begin{enumerate}[(i)] \item $(M_t^\eta\circ \pi)_{t\in[0,1]}$ is a continuous $(\mathcal{H}^\gamma_t)-${martingale}, on the complete probability space $(W\times S, 
 \mathcal{H}_1^\gamma, \gamma),$ where $\mathcal{H}_t:=( \mathcal{B}_{t+}^0(W)\otimes \mathcal{B}_t(S)),$ for all $t\in[0,1]$.
  \item For any (not necessarily right-continuous) bounded $(\mathcal{F}_t^\eta)-${martingale} $(M_t)$, on the probability space $(W,\mathcal{B}(W)^\eta, \eta)$, $(M_t\circ \pi)$ is both a $(\mathcal{H}_t^\gamma)-$martingale, and a $(\mathcal{H}^\gamma_{t+})-${martingale}, on the complete probability space $(W\times S, \mathcal{H}_1^\gamma,\gamma)$.
 \item  $\gamma \in \Pi_c(\eta, Y_\star \mathcal{P})$.
 \end{enumerate}
 \end{proposition}
\nproof
For the sake of clarity, we write the proof with $d=1$ ; recall that $\mathcal{B}_t(E)=\mathcal{B}_t(W)= \mathcal{B}_{t+}^0(W),$ so that $\mathcal{H}_t =\mathcal{B}_t(E)\otimes \mathcal{B}_t(S)$, for all $t\in[0,1]$. Moreover, the weak uniqueness of solutions to~(\ref{SDEJ}) implies $X_\star\mathcal{P}=\eta$.  We now start the proof. Since any right-continuous $(\mathcal{H}_t^\gamma)-${martingale} is a $(\mathcal{H}_{t+}^\gamma)-${martingale}, from the predicable representation property, and from the existence of continuous modifications to any $(\mathcal{F}_t^\eta)-$martingale,  we obtain that $(ii)$ follows from $(i)$. By Lemma~\ref{consistentreg}, since $\mathcal{B}_1(E)= \mathcal{B}(W)$, $(iii)$ follows from $(ii)$ ;  it is enough to prove $(i)$. Let $M : \Omega \to W$ and $u:\Omega \to W$ be  two $
\mathcal{A} / \mathcal{B}(W)-$measurable maps, such that $\mathcal{P}-a.s.$ $M:= \int_0^. \sigma_s(X)dB_s$ and $u:= \int_0^. b_s(X) ds$. By~(\ref{SDEJ}), $M$ is $(\mathcal{G}
_t^X)-$ adapted, where $\mathcal{G}_t^X:= X^{-1}(\mathcal{B}_{t}(W))$, for all $t\in [0,1]$. Hence, by classical results on transformations of measure (see Proposition 1.3. of \cite{CRUZLAS} for a statement, with the same notation),  there exists a $\mathcal{B}(W)^\eta / \mathcal{B}(W)-$measurable map $\widehat
{M} : W\to W$, such that $\mathcal{P}-a.s.,$ $M=\widehat{M}\circ X$. In particular the process $(\widehat{M}_t)$ associated to $\widehat{M}$ is a $(\mathcal{F}_t^\eta)-${martingale} 
on $(W, \mathcal{B}(W)^\eta, \eta)$, where $\widehat{M}_t:= W_t \circ \widehat{M}$ (pullback of maps). Similarly there exists some measurable map $\widehat{u} : W\to W$,  such that $\mathcal{P}-
a.s.$ $\widehat{u} \circ X := \int_0^. b_s(X) ds$, where $t\to \widehat{u}_t:= W_t\circ \widehat{u} \in \mathbb{R}$ is a $(\mathcal{F}_t^\eta)-$adapted process, on the probability space $(W,\mathcal{B}(W)^\eta, \eta)$, which is absolutely continuous outside a $\eta-$negligible set. Thus, from~(\ref{SDEJ}), we obtain  $$W_t= W_0 + M_t^\eta + \int_0^t b_t(\omega) ds = W_0+ \widehat{M}_t  + \widehat{u}_t, \ \forall t \ \  \eta-a.s.,$$ so that $M = \widehat{M}\circ X = M^\eta \circ X \ \mathcal{P}-a.s.,$ where we used that any continuous 
{martingale} of finite variation, starting at $0$, vanishes. On the other 
hand, since $Y$ is an $(\mathcal{A}_t)-$adapted map, and $X$ is $(\mathcal{A}_t)-$adapted, we obtain $\mathcal{H}^{X,Y}_t \subset \mathcal{A}_t,$ for all $t\in[0,1]$, where
$\mathcal{H}^{X,Y}_t:= \sigma(X^{-1}(\mathcal{B}_t(W)) \cup Y^{-1}(\mathcal{B}_t(S)))^\mathcal{P}.$  Whence,  we finally obtain $$E_{\gamma}[(M_t^\eta\circ \pi- M_s^\eta\circ \pi) 1_A\circ \pi 1_B\circ \widetilde{\pi}] = E_{\mathcal{P}}[(M_t-M_s) 1_{X^{-1}(A)}1_{Y^{-1}(B)}] = 0,$$ for all $A\in 
\mathcal{B}_s(W)$, $B\in \mathcal{B}_s(S),$ and $s\leq t$, since $(M_t)$ is a $(\mathcal{A}_t)-$martingale. By a monotone class argument, it yields $(i)$. \nqed
\begin{proposition}
 \label{propja2}
Under the same assumptions as Proposition~\ref{propja1},  we still denote by $\eta$ the unique law of solutions to~(\ref{SDEJ}).  Further assuming that $S$ is endowed with a regular filtration, for all $\nu\in \mathcal{P}_{S}$, and $\gamma \in \Pi_c(\eta,\nu)$, there exists a {complete stochastic basis} $(\Omega, \mathcal{A}, (\mathcal{A}_t)_{t\in[0,1]}, \mathcal{P})$, and 
\begin{enumerate}[(i)]
\item $(X_t)_{t\in[0,1]}$ (resp. $(B_t)_{t\in[0,1]}$) an $(\mathcal{A}_t)-$adapted continuous process (resp. an $(\mathcal{A}_t)-$Brownian motion), such that $(X, B)$ solves~(\ref{SDEJ}), on this {complete stochastic basis},
\item a measurable $(\mathcal{A}_t)-$adapted map $Y : \Omega \to S$, i.e. $Y^{-1}(\mathcal{B}_t(S)) \subset \mathcal{A}_t,$ for all $t\in[0,1]$,
\end{enumerate}
such that  $$\gamma = (X\times Y)_\star \mathcal{P},$$  $X: \Omega \to W$ denoting a $\mathcal{A}/\mathcal{B}(W)-$measurable map, associated to  $(X_t)$. Moreover, if $det
(\sigma_t)\neq 0$ for all $t\in[0,1]$, $\eta-a.s.$, then we can take $(\Omega, \mathcal{A},(\mathcal{A}_t)_{t\in[0,1]}, \mathcal{P})=(W\times S, (\mathcal{B}(W)\otimes \mathcal{B}(S))^\gamma, 
\mathcal{G}_{t+}^\gamma, \gamma)$, $X=\pi$, $Y=\widetilde{\pi}$, where $(\mathcal{G}_t)=(\mathcal{B}_t(W)\otimes \mathcal{B}_t(S))$. \end{proposition}
\nproof
As in Proposition~\ref{propja1}, for the sake of clarity, we write the proof for $d=1$. Define  $\widetilde{\Omega}:= W\times S$, $\widetilde{X}:= \pi$ (projection on the first component of $W\times S$), $\widetilde{Y}:= \widetilde{\pi}$ (projection on the second  component of $W\times S$), $\widetilde{\mathcal{A}} := (\mathcal{B}(W)\otimes \mathcal{B}(S))^\gamma$, $(\widetilde{\mathcal{G}}_t) := (\mathcal{B}_t(W)\otimes \mathcal{B}_t(S))$, $\widetilde{\mathcal{P}}:= \gamma$. Since~$\widetilde{X}_\star \widetilde{\mathcal{P}} = \eta$,~(\ref{coordMRT}) yields   \begin{equation} \label{25add} \widetilde{X}_t= \widetilde{X}_0 + M_t^{\widetilde{X}} + \int_0^t b_s(\widetilde{X})ds  \ for \ all \ t\in[0,1] \ \widetilde{\mathcal{P}}-a.s.,\end{equation} where $M^{\widetilde{X}}:= M^\eta\circ \widetilde{X}$ ;  in particular Dolean's approximations for the predicable covariation process ensure that $\mathcal{P}-a.s.$,  $<M^{\widetilde{X}}>= \int_0^. \sigma^2_s(\widetilde{X}) ds.$ Moreover, by $(iii)$ of Proposition~\ref{characsmn}, $$E_{\widetilde{\mathcal{P}}}[M_t^{\widetilde{X}}- M_s^{\widetilde{X}}| \widetilde{\mathcal{G}}^\gamma_s] = E_{\eta}[M_t^\eta- M_s^\eta|\mathcal{F}_s^\eta] \circ \widetilde{X} = 0 \ \widetilde{\mathcal{P}}-a.s.,$$ for all $s<t$. Thus, $(M_t^{\widetilde{X}})$ is a $(\widetilde{\mathcal{G}}_t^\gamma)-$martingale ; since it is continuous, it is also a $(\widetilde{\mathcal{A}}_t)-$martingale, where $(\widetilde{\mathcal{A}}_t):=(\widetilde{\mathcal{G}}^\gamma_{t+}).$ If, for all $s\in[0,1]$, $\sigma_s$ is invertible $\eta-a.s.$, then the result directly follows, from  Theorem 7.1. of \cite{I-W}, which is the elementary version of the martingale representation theorem. Otherwise, define $\Omega := \widetilde{\Omega} \times W$, $\mathcal{A}:= \widetilde{\mathcal{A}} \otimes \mathcal{B}(W)$, $\mathcal{P}:= \widetilde{\mathcal{P}}\otimes \mu$ (standard Wiener measure ; the law of standard Brownian motions), $(\mathcal{G}_t) := (\widetilde{\mathcal{A}}_t\otimes \mathcal{B}_{t}^0(W))$,  $(\mathcal{A}_t):= (\mathcal{G}_{t+}^{\mathcal{P}})$, $X:= \widetilde{X}\circ \pi^\Omega$, and $Y:= 
\widetilde{Y}\circ \pi^\Omega$, where $\pi^\Omega : (\omega, \widetilde{\omega}) \in \widetilde{\Omega} \times W \to \omega \in \widetilde{\Omega}$. Similarly, by~(\ref{25add}), the proof of the sharp {martingale representation theorem} (see p.89-91 of \cite{I-W}), ensures the existence of a $(\mathcal{A}_t)-${Brownian motion} $(B_t)$, such that the statement is satisfied.\nqed
\begin{remarkk}
\label{meprl}
Proposition~\ref{propja1} and Proposition~\ref{propja2}  apply together with Example~\ref{example1}.  Take $S= [0,1]$ (resp. $W$, resp. $[0,1]^n$), with the {regular 
filtrations} given respectively by Example~\ref{example1},  and denote by $\tau$ (resp. $Z$, resp. $\tau_1\times ...\times \tau_n$) the map $Y$ in the statements of Proposition~\ref{propja1} and Proposition~\ref{propja2}. Since  $(\Omega, \mathcal{A}, (\mathcal{A}_t), \mathcal{P})$ is a {complete stochastic basis}, by Example~\ref{example1}, the statement 
that $Y$ is an $(\mathcal{A}_t)-$adapted map reads : $\tau$ is an $(\mathcal{A}_t)-${stopping time} (resp. the stochastic process $(Z_t)$ is $(\mathcal{A}_t)-$adapted, resp. for all $i$, $\tau_i$ is a $(\mathcal{A}_t)-$stopping time). Subsequently, we will merely refer to Proposition~\ref{propja1} and Proposition~\ref{propja2}. \end{remarkk}

\section{Applications to stochastic control}\label{5} On suitable spaces, fixed in the next subsection, Corollary~\ref{prim} (resp. Corollary~\ref{prim2}) provides a straightforward equivalence of causal optimization problems, introduced in the first part of the paper, to \textit{optimal stopping problems} (resp. to \textit{extended optimal Skorokhod embedding problems} ; see \cite{BHC}). Our approach holds in \textit{non-Markovian} frameworks, and  the equivalence is directly obtained by projections on component of $E\times S$. Thus, concerning optimal \textit{SEP}, this approach, which uses a straightforward and different representation than \cite{BHC}, where they use auxiliary stopping times, improves known results on the primal attainment for those problems ; it also clarifies the role of objects called $RST$ in \cite{BHC}. However, we do not consider the geometry of optimal plans, which is the main part of the latter ; it would be interesting to investigate, whether the short approach of the present paper, could also enlighten these geometric aspects.  

\subsection{Framework}
In the whole section, $E=W$ (resp. $S=[0,1]$), $I=[0,1]$,  $(\mathcal{B}_t(E))$ (resp. $(\mathcal{B}_t(S))$ is given by~(\ref{defilW}) (resp. by $(\mathcal{B}_t([0,1])$) given in Example~\ref{example1}) ; see subsection~\ref{soussection42}.
\subsection{Equivalence to optimal stopping problems} Recall that $\mathcal{P}_c^{[0,1]}(\eta):= \cup_{\nu\in \mathcal{P}_{[0,1]}} \Pi_c(\eta,\nu)$ ; $(\rho_t^W)$ is given by~(\ref{rhotWdefcor}).
\begin{corollary}
 \label{prim}
Under the assumptions of Proposition~\ref{propja2} on~(\ref{SDEJ}), still denote by $\eta$ the law of its solutions. Given a non-negative l.s.c. map $c_a : W\times [0,1] \to \mathbb{R} \cup\{+\infty\},$ define $c:= c_a\circ \phi$, the pullback of the map $c_a$, with the map   $\phi : (\omega, t) \in W\times [0,1] \to (\omega_{.\wedge t}, t) \in W\times [0,1].$ The primal problem  \begin{equation} \label{563prov} P_{\eta}:=\inf\left(\left\{ \int_{W\times [0,1]} c(\omega,t) d\gamma(\omega,t) \middle| \gamma \in \mathcal{P}_c^{[0,1]}(\eta) \right\}\right),\end{equation} is equivalent to minimize \begin{equation} \label{primalOS}J_\eta:=  \inf(E_{\mathcal{P}}[c_a(\rho^W_{\tau}(X),\tau)]),\end{equation} where the infimum of~(\ref{primalOS}) is taken on all the {complete stochastic basis} $(\Omega, \mathcal{A}, (\mathcal{A}_t)_{t\in[0,1]}, \mathcal{P})$, such that there exists an $(\mathcal{A}_t)-${Brownian motion} $(B_t)$, an $(\mathcal{A}_t)-${stopping time} $\tau$ (bounded by $1$), and an $(\mathcal{A}_t)-$ adapted continuous process $(X_t)$, which satisfy \begin{equation} \label{stardercor} X_t = X_0 + \int_0^t \sigma_s(X) dB_s + \int_0^t b_s(X) ds ; Law(X_0) = \eta_0 \ , \ for \ all \ t\in[0,1] \ \mathcal{P}-a.s..\end{equation}That is, $J_\eta= P_\eta$, and $\gamma\in \mathcal{P}_c^{[0,1]}(\eta)$ attains the infimum of~(\ref{563prov}), if and only if, there exists a pair $(X,\tau)$, on a complete stochastic basis $(\Omega, \mathcal{A}$, $(\mathcal{A}_t)_{t\in[0,1]},\mathcal{P})$, which attains~(\ref{primalOS}), for some $(\mathcal{A}_t)-$Brownian motion $(B_t)$, and satisfies $\gamma = (X\times \tau)_\star \mathcal{P}$. Moreover, the infimum of~(\ref{primalOS}) is attained.
  \end{corollary}
  \nproof
  Since $\phi$ is continuous,  $c:= c_a \circ \phi,$ is non-negative and lower semi-continuous. Together with Remark~\ref{meprl}, we obtain $P_\eta\leq J_\eta$ (resp. $J_\eta \leq  P_\eta$) by Proposition~\ref{propja1} (resp. by Proposition~\ref{propja2}), so that $P_\eta =J_\eta$, and the optimum are identified as stated.  Whence, the attainment of~(\ref{primalOS}), follows from Corollary~\ref{ToSEP}. 
  \nqed
 \subsection{An alternative to the Beiglboeck-Cox-Huesmann representation of  extended SEP}
    \begin{corollary}
 \label{prim2}
Under the assumptions and notations of Corollary~\ref{prim}, we further denote by $e$ the evaluation map $e : (\omega, t) \in W\times [0,1] 
 \to e(\omega, t) :=\omega(t) \in \mathbb{R}^d$. Given a non-negative l.s.c. map  $c_a : W\times [0,1] \to \mathbb{R} \cup\{+\infty\},$ we still associate the map $c$ of Corollary~\ref{prim}.  Then, in the same acceptation as Corollary~\ref{prim}, for all $\widetilde{\mu} \in \mathcal{P}_{\mathbb{R}^d}$, the primal problem $$P^{e}_{\eta,\widetilde{\mu}}:=  \inf\left(\left\{\int_{W\times [0,1]}c(\omega,t) d\gamma (\omega,t) \middle| \gamma \in \mathcal{P}^{[0,1]}_c(\eta),  e_\star \gamma = \widetilde{\mu}\right\}\right),$$  is equivalent  to  \begin{equation} \label{primSEPldmf}J_{\eta,\widetilde{\mu}}:= \inf(E_{\mathcal
 {P}}[c_a(\rho^W_{\tau}(X),\tau)] :  Law(X_{\tau}) = \widetilde{\mu}), \end{equation} where the infimum is taken on all the {complete stochastic basis} $(\Omega, \mathcal{A}, (\mathcal{A}_t)_{t\in
 [0,1]}, \mathcal{P})$, such that there exists an $(\mathcal{A}_t)-${Brownian motion} $(B_t)_{t\in[0,1]}$, an $(\mathcal{A}_t)-${stopping time} $\tau$ (bounded by $1$), and an $
 (\mathcal{A}_t)-$adapted continuous process $(X_t)_{t\in[0,1]}$, which meet~(\ref{stardercor}), and $$(X_\tau)_\star \mathcal{P} = \widetilde{\mu}.$$ Moreover, if $P_{\eta,\widetilde{\mu}}< \infty,$ then the infimum of~(\ref{primSEPldmf}) is attained.
  \end{corollary}
 \nproof
Since $e$ is continuous, the result follows from  Corollary~\ref{ToSEP}, similarly to the proof of Corollary~\ref{prim}.
\nqed
\section{Stochastic differential equations as causal optimal transportation problems}
 \label{6} This section investigates connections to \textit{stochastic differential equations}. The first subsection fixes the notation for the whole section ; Proposition~\ref{Rbscharac}, and Proposition~\ref{prop423}, formulate  results on transformations of the law of  standard {Brownian motions}, denoted by $\mu$, in the framework of the first part of the paper. This yields the equivalence between certain \textit{stochastic differential equations} and \textit{causal Monge-Kantorovich problems}.  Proposition~\ref{Rbscharac} investigates some deterministic transport, from laws of continuous stochastic processes to $\mu$, under the symmetric constraint of Section~\ref{4}. It provides a concise proof to Proposition~\ref{prop423}, which is a synthesis of several well known results on transformations of $\mu$, within the analytic framework of Section~\ref{1}-Section~\ref{3}. Lemma~\ref{lemmaweak} provides an analytic reformulation of  representation formulas of the \textit{entropy} with respect to $\mu$. It is a complement to the celebrated formula of   \cite{F1} ; the latter has deep implications in so-called the \textit{Schr\"{o}dinger problems} (see  \cite{leo}, \cite{F1}, \cite{Zamb1}). In \cite{ASU-3}, \cite{ASULAST}, a strong version of this formula was used to investigate strong existence problems for stochastic differential equations ; it is closely related to \cite{DBF} (see \cite{ASULAST}) ; it was weakened in  \cite{LEHEC},\cite{FISHER}, \cite{RLASU},  \cite{ABSAPI}. Apart the analytic formulation, the novelty here is the restriction of the optimization to the symmetric constraint of Section~\ref{4}. Lemma~\ref{lemmaweak} provides a reformulation of these results, and investigates the related optimums in connection with \textit{Malliavin calculus}.  Proposition~\ref{Talag} shows that the previous Lemma~\ref{lemmaweak}  trivially entails some \textit{Talagrand's inequality} (see \cite{ASU-3}, \cite{Tal}).  In this particular case, this suggests to see  the difference between the \textit{relative entropy} and the \textit{Wasserstein distance}, as the price to pay to buy all the information contained in $\mathcal{B}(W)$, at $t=0$. Although the proof appeared somewhere else (in \cite{RLASU}, \cite{ASU-3}), it was not written explicitly in the form of optimal transport until now.   Finally, from Lemma~\ref{lemmaweak}, Theorem~\ref{SDEMK} characterizes solutions to certain stochastic differential equations as optimum to \textit{causal Monge-Kantorovich problems} ;  the existence of a unique strong solution is related to the corresponding \textit{Monge problem}.
\subsection{Framework of this section} \textit{Until the end of the paper, within the framework of the first part of the paper, we take $E= W$, $S= W$, $I=[0,1]$,  and $(\mathcal{B}_t(E))$ (resp. $(\mathcal{B}_t(S))$ is given by~(\ref{defilW}) (resp. by the natural filtration~(\ref{evalprokf}));}  $(\mathcal{F}_t^\eta)$ is given by~(\ref{emprol}). Thus, we are in the specific framework of subsection~(\ref{symsubsec}), and given $\eta,\nu\in \mathcal{P}_W$, $\mathcal{R}_{as}(\eta,\nu)$ is defined by~(\ref{raz}), with $Z=W$. Given $U\in \mathcal{R}(\eta,\nu)$, and any $\mathcal{B}(W)^\eta /\mathcal{B}(W)$-measurable map $\widetilde{U}: W\to W$, whose related $\eta-$equivalence class of maps is $U$, set $$U_t : \omega \in W \to U_t(\omega):=W_t(\widetilde{U}(\omega)) \in \mathbb{R}^d,$$ where $W_t : \omega \in W \to W_t(\omega):=\omega(t) \in \mathbb{R}^d$, for all $t\in [0,1]$. Subsequently, we refer to such continuous process $(U_t)$, on the probability space $(W,\mathcal{B}(W)^\eta,\eta)$, as a \textit{process associated to $U$}. By Example~\ref{example1}, and Proposition~\ref{consistentreg2},  $U\in \mathcal{R}_a(\eta,\nu)$, if and only if,  any (and then all) process 
$(U_t)$ associated to $U$, is $(\mathcal{F}_t^\eta)$-adapted. Thus, in this particular case, \textit{causal transference plans can be interpreted, as a relaxation of adapted (also called causal) processes} on canonical spaces ; whence they get their name.  We consider transformations of the {Wiener measure}, which we denote by $\mu\in \mathcal{P}_W$ ; the \textit{law of standard Brownian motions}, seen as random continuous paths. That is, the unique element of $\mathcal{P}_W$, such that the evaluation process $(W_t)$ is a {standard Brownian motion} on $(W,\mathcal{B}(W), \mu)$. Since it is a continuous process, $(W_t)$ is also a $(\mathcal{F}_t^\mu)-${Brownian motion} on $(W,\mathcal{B}(W)^\mu,\mu)$.
Finally, recall the \textit{Sobolev derivative} extends to the infinite dimensional space $W$, endowed with the measure $\mu$ ; to handle Banach spaces, equivalence classes of maps must be preserved by translations, which entails \textit{quasi-invariance} issues. Thus, denote by $H$, the so-called \textit{Cameron-Martin space} of $\widetilde{\omega}\in W$ such that ${\tau_{\widetilde{\omega}}}_\star\mu \sim \mu$, \textit{i.e.} equivalent,  where $\tau_{\widetilde{\omega}} : \omega \in W \to \omega + \widetilde{\omega} \in W$. Due to the Cameron-Martin theorem, it is actually given by  $$H= \left\{ h\in W,  h= 
\int_0^. \dot{h}_s ds \middle| \int_0^1 |\dot{h}_s|_{\mathbb{R}^d}^2 ds <\infty  \right\},$$  which is an Hilbert space for the scalar product $$<h,k>_H :=\int_0^1 <\dot{h}_s,\dot{k}
_s>_{\mathbb{R}^d} ds,$$ for $h,k\in H$. As a classical application of the It\^{o}-Nisio theorem (see \cite{ITN}), \textit{ $|.|_H$ extends to a non-negative l.s.c. map on $W$}, by setting 
\begin{equation} |.|_H : \omega \in W \to | \omega|_H= \begin{cases}
        \sqrt{<\omega,\omega>_H} & \text{if } \ \omega \in H \\
       + \infty & \text{otherwise}
       \end{cases} .\end{equation} 
 Below, \textit{causal Monge Kantorovich problems} are investigated, for a cost map $$c :(x,y)\in W\times W \to c(x,y):= |x-y|_H^2 \in \mathbb{R}\cup \{+\infty\}.$$ \textit{A reader not familiar with Malliavin calculus is encouraged to skip the end of this subsection}. We refer to \cite{MAL} and  \cite{I-W}, for an overview on this topic. Recall this derivative is first defined, as the $H-derivative$ (see \cite{KUO}), on the set of smooth polynomials $\mathcal{P}_{ol}$ on $W$ (see \cite{I-W}), so that to obtain a closable operator $\nabla : \mathcal{P}_{ol}\subset L^2(\mu) \to L^2(\mu,H).$ Denoting $\DD_{2,1}$ the completion of $\mathcal{P}_{ol}$ with respect to the norm of the graph $||.||_{2,1} : F \in \mathcal{P}_{ol}  \to ||F||_{2,1} := |F|_{L^2(\mu)} + |\nabla F|_{L^2(\mu,H)}\in \mathbb{R},$ as a closable operator, $\nabla$ extends to a map $$ \nabla : F\in \DD_{2,1}\subset L^2(\mu) \to \nabla F \in L^2(\mu,H) ;$$ the so-called \textit{Malliavin derivative}. Take $F\in L^2(\mu)$, which is the $\mu-$ equivalence class of maps of some  $f \in \mathcal{L}^2(\mu).$  Assuming that $f$ is $H-differentiable$ (see \cite{KUO}), and that its derivative defines an element of $L^2(\mu,H)$, we have  $$<\nabla F, h>_H = \frac{d}{d\lambda}f(\omega + \lambda h)|_{\lambda =0} \ \mu-a.s.,$$ for all $h\in H$.  For any $X\in \DD_{2,1}$, we denote by $(D_s X)$ the derivative of $\nabla X$, with respect to the Lebesgue measure ; it satisfies $$ \nabla X = \int_0^. D_s X ds \ \mu-a.s..$$ The density of $\mathcal{P}_{ol}$ in $L^2(\mu)$ ensures the existence of the adjoint to the Malliavin derivative $\nabla$, the so-called \textit{divergence}, which coincides with the stochastic integral on the subset of  elements of its domain (subset of $L^2(\mu,H)$), whose related processes are $(\mathcal{F}_t^\mu)-$adapted. Together with the \textit{martingale representation theorem}, it yields the so-called \textit{Clark-Ocone formula}, which reads $$X= E_\mu[X] + \int_0^1 E_\mu\left[D_sX \middle| \mathcal{F}_s^\mu\right]dW_s. \ \mu-a.s.,$$ for all $X\in \DD_{2,1}$. Subsequently, for convenience of notations, given  $X,Y\in M_{\mathcal{P}}((\Omega, \mathcal{A}), (W, \mathcal{B}(W))$, we denote by $(X\times Y)_\star\mathcal{P}$, the joint law $(\widetilde{X}\times \widetilde{Y})_\star \mathcal{P}$ of any (and then all), pair of measurable maps $\widetilde{X}, \widetilde{Y} : \Omega \to W$, whose $\mathcal{P}-$equivalence class is $X$ (resp. $Y$). Finally, for $X\in M_{\mathcal{P}}((\Omega, \mathcal{A}), (W, \mathcal{B}(W))$, and $\widetilde{X} : \Omega \to W$, by $\mathcal{P}-a.s.$ $X= \widetilde{X}$, we denote that $X$ is the $\mathcal{P}-$equivalence class of $\widetilde{X}$.

\subsection{Transformations of the law of standard Brownian motions}\label{sub3SDE} In the statements below, given $\nu\in \mathcal{P}_W$, $\mathcal{R}_{as}(\nu,\mu)$ (resp. $\Pi_{cs}(\mu,\nu)$), is given by~(\ref{raz}) (resp. by ~(\ref{causcoupl})), of Section~\ref{4}.
   \begin{proposition}
 \label{Rbscharac}
 Given $\nu \in \mathcal{P}_W$ and $V\in \mathcal{R}(\nu,\mu)$, we have $V\in \mathcal{R}_{as}(\nu,\mu)$ if and only if $(V_t)$ is a $(\mathcal{F}_t^\nu)-${Brownian motion} on the probability space $(W, \mathcal{B}(W)^\nu,\nu)$, for any (and then all) continuous process $(V_t)$ associated to $V$.
 \end{proposition}
 \nproof
 From the definition, $V\in \mathcal{R}_{as}(\nu,\mu)$, if and only if, both\begin{equation} \label{condone} j(V):= (I_W\times V)_\star \nu \in \Pi_c(\nu,\mu)\end{equation} and \begin{equation} \label{condone2} R_\star j(V):= (V \times I_W)_\star \nu \in \Pi_c(\mu,\nu), \end{equation} hold,  $j$ (resp. $R$)  denoting the map given by~(\ref{jdef}) (resp. by~(\ref{Rdef})).  By Proposition~\ref{consistentreg2},~(\ref{condone}) is equivalent to $$\widetilde{V}^{-1}(\mathcal{B}^0_t(W))\subset  \mathcal{B}_{t+}^0(W)^\nu = \mathcal{F}_t^\nu$$ for all $t\in[0,1]$,  $\widetilde{V}: W\to W$ denoting any $\mathcal{B}(W)^\nu / \mathcal{B}(W)-$measurable map, whose $\nu-$equivalence class is $V$.  Since $\mathcal{R}_{as}(\nu,\mu):=j^{-1}(\Pi_{cs}(\nu, \mu))$, by Proposition~\ref{propja2} and Proposition~\ref{propja1}, with $\sigma_t= I_{\mathbb{R}^d}$ and $b_t=0$, the second equality is equivalent to $(V_t)$ is a $(\mathcal{G}_{t+}^\nu)$-{Brownian motion} on $(W,\mathcal{B}(W)^\nu,\nu)$, $\mathcal{G}_t$ denoting the sigma-field $\sigma(\widetilde{V}^{-1}(\mathcal{B}_{t+}^0(W)) \cup \mathcal{B}_{t}^0(W )),$ for all $t\in [0,1]$. Since $(\mathcal{F}_t^\nu)$ is right-continuous, the result follows.
  \nqed
\begin{proposition}
 \label{prop423}
  For any $\nu\in \mathcal{P}_W$, absolutely continuous probability, with respect to the Wiener measure $\mu$, there exists a unique $V^\nu\in \mathcal{R}_{as}(\nu,\mu)$,  such that  \begin{equation} \label{vnudef} V^\nu = I_W - \int_0^. v_s^\nu ds \ \nu-a.s., \end{equation} where $(v_s^\nu)$ is a $(\mathcal{F}_t^\nu)-$predicable process, which satisfies \begin{equation} \label{girsint} \int_0^1 |v_s^\nu|_{\mathbb{R}^d}^2 ds <\infty \ \nu-a.s., \end{equation}  and where $I_W :\omega \in W \to \omega \in W$ denotes the identity map on $W$. Moreover the following hold,
  \begin{enumerate}[(i)]
  \item  We have   \begin{equation}  \label{folm} 2\mathcal{H}(\nu | \mu) = \int_{W\times W} |x-y|_H^2 d\gamma_c(x,y),\end{equation}  where $\gamma_c:= (V^\nu \times I_W)_\star \nu,$ and where \begin{equation} \label{entrodefini} \mathcal{H}(\nu| \mu):=E_\nu\left[\ln \frac{d\nu}{d\mu}\right]\end{equation} denotes the relative entropy. Moreover, $\gamma_c \in \Pi_{cs}(\mu,\nu) \subset \Pi_c(\mu,\nu)$.    \item $(W_t,V_t^\nu)$ is a {weak solution} to $$dX_t = dB_t + v_t^\nu\circ X dt ; X_0=0,$$ on the probability space $(W,\mathcal{B}(W)^\nu,\nu)$, with the filtration $(\mathcal{F}_t^\nu)$, for any (and then all) process $(V_t^\nu)$ associated to $V^\nu$. 
  \item the Radon-Nikodym derivative of $\nu$ is given by  \begin{equation} \label{densitygirs} \frac{d\nu}{d\mu} = \exp\left(\int_0^1 v_s^\nu dW_s  -\frac{1}{2}\int_0^1|v_s^\nu|_{\mathbb{R}^d}^2ds\right) \ \nu-a.s. \end{equation}
  \item further assuming that $\frac{d\nu}{d\mu}\in \DD_{2,1},$ \begin{equation}  \label{empo} v_t^\nu = E_\nu\left[D_t \ln \frac{d\nu}{d\mu} \middle| \mathcal{F}_t^\nu \right] \   dt\otimes d\nu-a.s..\end{equation} 
     \end{enumerate}
     We call $V^\nu$ the Girsanov shift of $\nu$. 
\end{proposition} 
 \nproof
 Together with the {martingale representation theorem}, the {Girsanov theorem}  ensures the existence of a unique $(\mathcal{F}_t^\nu)-${Brownian motion} $(V_t^\nu)$, on the probability space $(W,\mathcal{B}(W)^\nu,\nu)$, such that $$W_t = V^\nu_t + \int_0^t v_s^\nu ds,$$ for  all  $t\in [0,1]$, $\nu-a.s.,$ where $(v_s^\nu)$ is a $(\mathcal{F}_t^\nu)-$predicable process which satisfies~(\ref{girsint}), and such that~(\ref{densitygirs}) holds (see \cite{F1} and the references therein) ; in particular it implies $(ii)$. Thus, in (\ref{densitygirs}), the stochastic integral is with respect to a {semi-martingale}.  From Proposition~\ref{Rbscharac}, we obtain $V^\nu\in \mathcal{R}_{as}(\nu,\mu)$, and the uniqueness follows from the fact that a {continuous martingale} of finite variation vanishes. By the celebrated {representation formula of the entropy} of \cite{F1}, we obtain \begin{equation} \label{ABSFOLLMER} 2\mathcal{H}(\nu | \mu)= E_\nu\left[\int_0^1 |v_s^\nu|_{\mathbb{R}^d}^2 ds \right] = E_\nu[|V^\nu -I_W|_H^2] , \end{equation}  $|.|_H$ denoting the Cameron-Martin norm ; this proves~(\ref{folm}). Assuming that $\frac{d\nu}{d\mu} \in \DD_{2,1}$, the Clark-Ocone formula yields  $$\frac{d\nu}{d\mu} =  1 + \int_0^1 E_\mu\left[D_s \frac{d\nu}{d\mu}  \middle| \mathcal{F}_t^\mu \right] dW_s  \ \mu-a.s.,$$ and thus $\nu-a.s.$. Whence  \begin{equation} \label{densm} \frac{d\nu}{d\mu} = \exp\left( \int_0^1 E_\nu\left[D_s \ln \frac{d\nu}{d\mu} \middle| \mathcal{F}_s^\nu \right] dW_s - \frac{1}{2} \int_0^1 \left|E_\nu\left[D_s \ln \frac{d\nu}{d\mu} \middle| \mathcal{F}_s^\nu \right] \right|_{\mathbb{R}^d}^2 ds \right)  \ \nu-a.s., \end{equation} stemming from It\^{o}'s formula, together with the condition $\frac{d\nu}{d\mu}>0 \ \nu-a.s.$, and with the fact that $\nabla$ is a local operator. Since any {martingale}, starting from $0$,  of finite variations, vanishes,~(\ref{empo}) follows from~(\ref{densitygirs}) and~(\ref{densm}). Finally, by Proposition~\ref{Rbscharac}, $V^\nu\in \mathcal{R}_{as}(\nu,\mu)= j^{-1}(\Pi_{cs}(\nu,\mu))$. Thus, by symmetry $\gamma_c :=R_\star j(V^\nu) \in \Pi_{cs}(\mu,\nu)$. \nqed

 \subsection{Stochastic differential equations as optimal transport problems} For $\nu \in \mathcal{P}_W$, we define its relative entropy $\mathcal{H}(\nu | \mu)$, with respect to the Wiener measure $\mu$, by~(\ref{entrodefini}) if $\nu<<\mu$ (i.e. absolutely continuous), and by $\mathcal{H}(\nu | \mu) = +\infty$ otherwise.
 \begin{lemma}
 \label{lemmaweak}
 For all $\nu\in \mathcal{P}_W,$ we have   \begin{equation} \label{entroweak} 2 \mathcal{H}(\nu| \mu) = \inf\left(\left\{ \int_{W\times W}|x-y|_H^2 d\gamma(x,y) \middle| \gamma \in \Pi_c(\mu,\nu) \right\}\right), \end{equation} $\mathcal{H}(\nu | \mu)$ denoting the relative entropy, and the infimum is attained by some $\gamma_{c} \in \Pi_{cs}(\mu,\nu) \subset \Pi_c(\mu,\nu)$. Moreover, if $\mathcal{H}(\nu | \mu)$ is finite, then $\gamma_{c}$ is unique, and it is given by \begin{equation} \label{optV} \gamma_{c}=  (V^\nu\times I_W)_\star\nu,\end{equation}  $V^\nu$ denoting the Girsanov shift of $\nu$ ; see Proposition~\ref{prop423}. In particular,~(\ref{entroweak}), still holds by substituting  $\Pi_{cs}(\mu,\nu)$ for $\Pi_{c}(\mu,\nu)$.
 \end{lemma}
 \nproof Denoting  $P_{\mu,\nu}$ the right hand term of~(\ref{entroweak}), Proposition~\ref
 {prop423} yields $P_{\mu,\nu} \leq 2 \mathcal{H}(\nu | \mu)$. Henceforth, we assume that $P_{\mu,\nu}$ is finite, and we take $\gamma\in \Pi_c(\mu,\nu)$, such that $E_\gamma[|\pi-\widetilde{\pi}|_H^2] <\infty.$ Define $u:=\int_0^. \dot{u}_s ds\in L^2(\gamma,H)$ by \begin{equation} \label{udefmp} \widetilde{\pi} = \pi + u \ \gamma-a.s..\end{equation} Since $\gamma=(\pi\times \widetilde{\pi})_\star \gamma \in \Pi_c(\mu,\nu)$, by Proposition~\ref{propja2}, $t\to W_t\circ \pi$ (resp. $t\to u_t$) is $(\mathcal{G}_{t+}^
 \gamma)-${Brownian motion} (resp. is a square integrable $(\mathcal{G}_{t+}^\gamma)-$adapted process), on the complete stochastic basis $(W\times W, \mathcal{B}(W\times W), (\mathcal{G}_{t+}^\gamma))$, where $(\mathcal{G}_t):= (\mathcal{B}_t(W)\otimes \mathcal{B}_t^0(W))$. By the Girsanov theorem, we obtain $\nu<<\mu$. On the other hand, by standard results on transformations of laws of {semi-martingales} (see Proposition 3.1. of \cite{CRUZLAS}, for a statement with the same notations), since $\widetilde{\pi}_\star \gamma = \nu$, we obtain  \begin{equation} \label{vcp} \int_0^. v_t^\nu\circ \widetilde{\pi} dt =\int_0^. E_\gamma[\dot{u}_t | \mathcal{G}_t^{\widetilde{\pi}}] dt \ \gamma-a.s.,\end
 {equation} where the right hand term denotes the dual predicable projection of $(u_t)$, on the $\gamma-$usual augmentation, of the natural filtration generated by $t\to W_t\circ \widetilde{\pi}$ ; see~(\ref{vnudef}). By Jensen's inequality, together with~(\ref{ABSFOLLMER}),  \begin{equation} \label{entroeqjens} 2\mathcal{H}(\nu | \mu) = E_{\gamma}\left[\int_0^1|E_\gamma[\dot{u}_t | \mathcal
 {G}_t^{\widetilde{\pi}}]|_{\mathbb{R}^d}^2 dt \right] \leq  E_{\gamma}\left[|\widetilde{\pi}- \pi|_H^2\right] \end{equation} follows from~(\ref{udefmp}) and ~(\ref{vcp}). Since in this case, the dual predicable projection is the orthogonal projection, on the closed linear subspace of adapted elements of $L^2(\mu,H)$,  the equality occur in~(\ref{entroeqjens}), if and only if,  $$\int_0^. v_t^\nu\circ \widetilde{\pi}=  \int_0^. \dot{u}_s ds \ \gamma-a.s..$$ By~(\ref{udefmp}) and~(\ref{vnudef}), since $\gamma=(\pi\times \widetilde{\pi})_\star \gamma$, the latter is equivalent to $\gamma = \gamma_c$, where $\gamma_c$ is given by~(\ref{optV}) ; by Proposition~\ref{prop423}, $\gamma_c  \in \Pi_{cs}(\mu,\nu) \subset \Pi_c(\mu,\nu)$. Thus,~(\ref{entroeqjens}) yields~(\ref{entroweak}), for all $\nu \in \mathcal{P}_W$. Moreover, assuming the entropy is finite, from the proof, the infimum is attained by $\gamma_c$, given by~(\ref{optV}). \nqed 

  \begin{proposition}
\label{Talag}
For any $\nu\in \mathcal{P}_W$, we have \begin{equation} \label{talinv}   \inf\left(\left\{ \int_{W\times W} |x-
y|_H^2 d\gamma(x,y) , \gamma \in \Pi(\mu, \nu)\right\}\right) \leq 2 \mathcal{H}(\nu | \mu). \end{equation}  \end{proposition}
\nproof
By Lemma~\ref{lemmaweak}, $\Pi_c(\mu,\nu) \subset \Pi(\mu,\nu)$ implies~(\ref{talinv}). \nqed

Theorem~\ref{SDEMK} provides an optimal transportations view, on strong existence of solutions to stochastic differential equations (see \cite{I-W}), which emphasizes the importance of the marginals, rather than the drift. Henceforth,  $v: (\omega,t) \in W\times [0,1] \to v_t(\omega) \in \mathbb{R}^d$ denotes  a bounded predicable map (see \cite{I-W}). To state Theorem~\ref{SDEMK}  below, recall some well known basic facts (for further details see \cite{I-W}).  As an application of the Girsanov theorem, the so-called \textit{transformation of the drift method} entails the weak existence and uniqueness of a solution to  \begin{equation} \label{SDEfin} dX_t = dB_t + v_t(X)dt ; X_0=0. \end{equation}  Denoting by $\nu$ the unique law of solutions to~(\ref{SDEfin}), it is the probability equivalent to $\mu$ such that  \begin{equation} \label{sdedens} \frac{d\nu}{d\mu} = \exp\left(\int_0^1 v_s dW_s -\frac{1}{2}\int_0^1 |v_s|_{\mathbb{R}^d}^2 ds\right) \ \mu-a.s. \end{equation}   Given $X,B :\Omega \to W$, two measurable maps, on a {complete stochastic basis} $(\Omega, \mathcal{A},(\mathcal{A}_t), \mathcal{P})$, let  $(X_t)$ (resp. $(B_t)$), be a pair of processes associated to $X$ (resp. to $B$)  which  is $(\mathcal{A}_t)-$adapted (resp. a $(\mathcal{A}_t)-${Brownian motion}).  From the very definition of a \textit{weak solutions}, $(X_t,B_t)$ solves~(\ref{SDEfin}) if and only if \begin{equation}  \label{BTX} (B\times X)_\star \mathcal{P}= (\widetilde{V}\times I_W)_\star \nu,\end{equation} $I_W$ denoting the identity map on $W$, and  $\widetilde{V}$ denoting the measurable map \begin{equation} \label{vtildef} \widetilde{V} : \omega \in W \to \omega-\int_0^. v_s(\omega) ds \in W.\end{equation} For convenience of notations we set $\gamma_\star := (\widetilde{V}\times I_W)_\star \nu$ ; we call it \textit{the unique joint law} of solutions to~(\ref{SDEfin}). In particular, $(\widetilde{V}_t)$ is a $(\mathcal{F}_t^\nu)-$Brownian motion, where $\widetilde{V}_t:= W_t\circ \widetilde{V}$, for all $t\in[0,1]$, $\circ$ denoting the pullback of maps. Finally, we say that~(\ref{SDEfin}) has a unique strong solution, if there exists a measurable map $F : W\to W$, such that, for all \textit{weak solution} $(B_t,X_t)$ to~(\ref{SDEfin}), on some {complete stochastic basis} $(\Omega,\mathcal{A}, (\mathcal{A}_t)_{t\in[0,1]}, \mathcal{P})$, we have $X= F(B) \ \mathcal{P}-a.s.,$ for the related maps and if, for all $(\mathcal{A}_t)-$Brownian motion $(B_t)$, on a {complete stochastic basis} $(\Omega,\mathcal{A}, (\mathcal{A}_t)_{t\in[0,1]}, \mathcal{P})$, the pair $(F(B),B)$ satisfies~(\ref{SDEfin}). 
 \begin{theorem}
 \label{SDEMK}
  Denoting by $\nu$  the probability given by~(\ref{sdedens}), the primal attainment of the  {causal Monge-Kantorovich problem}
   \begin{equation} \label{pmunu} P_{\mu,\nu} := \inf\left(\left\{ \int_{W\times W}|x-y|_H^2 d\gamma(x,y) \middle| \gamma \in \Pi_c(\mu,\nu) \right\}\right) \end{equation} is achieved  by a unique {causal transference plan}, which is $\gamma_{\star}$, the unique joint law of solutions to~(\ref{SDEfin}) ; the latter has a unique strong solution, if and  only if the optimal plan to~(\ref{pmunu}) induces a solution to the {causal Monge problem} \begin{equation} \label{pmonge} P^{Monge}_{\mu,\nu}= \inf\left(\int_W |x-U(x)|_H^2 d\mu(x), U\in \mathcal{R}_a(\mu,\nu)\right). \end{equation} Moreover, the above statement holds with $\Pi_{cs}(\mu,\nu)$ (resp. $\mathcal{R}_{as}(\mu,\nu)$), instead of  $\Pi_{c}(\mu,\nu)$ (resp. $\mathcal{R}_{a}(\mu,\nu)).$ Finally, by further assuming that $\frac{d\nu}{d\mu}\in \DD_{2,1}$, we have  $$\gamma_{\star}= ((I_W-\pi^\nu \nabla \ln \frac{d\nu}{d\mu})\times I_W)_\star\nu,$$ $\pi^\nu \nabla \ln \frac{d\nu}{d\mu}$ denoting the projection of $\nabla \ln \frac{d\nu}{d\mu}$, on the closed subspace of $L^2(\nu,H)$, whose elements are $(\mathcal{F}_t^\nu)-$adapted.
 \end{theorem}
 \nproof  Since a {continuous martingale} of finite variations vanishes, and $\nu$ is equivalent to $\mu$, $(iii)$ of Proposition~\ref{prop423} and~(\ref{sdedens}) yield $\nu-a.s.$ $V^\nu = \widetilde{V}$, $\widetilde{V}$ denoting~(\ref{vtildef}), and $V^\nu$ denoting the Girsanov shift of $\nu$, see Propostion~\ref{prop423}. In particular $\gamma_\star =\gamma_c,$ where $\gamma_c$ is defined by Proposition~\ref{prop423}. By~(\ref{ABSFOLLMER}), the entropy of $\nu$ w.r.t. $\mu$ is finite. Whence, it follows from Lemma~\ref{lemmaweak}, that~$\gamma_c \in  \Pi_{cs}(\mu,\nu)\subset \Pi_c(\mu,\nu)$, is the unique {causal transference plan} which attains $P_{\mu,\nu}$.   Finally, assuming the existence of a unique strong solution $F$,  since $(W_t)$ is a {Brownian motion} on $(W,\mathcal{B}(W)^\mu,\mu)$,  from the definition of the unique strong solution, $(F,I_W)$ is solution to~(\ref{SDEfin}) on $(W,\mathcal{B}(W)^\mu,\mu)$. Therefore,~(\ref{BTX}) implies $$\gamma_c = (I_W\times F)_\star \mu =j(U), $$  $U\in \mathcal{R}(\mu,\nu)$ denoting the $\mu-$equivalence class of maps associated to $F$, $j$ denoting the map~(\ref{jdef}). We obtain $U\in \mathcal{R}_a(\mu,\nu)$, from  $\gamma \in \Pi_c(\mu,\nu)$. Since $\gamma_c$ attains $P_{\mu,\nu}$, $U$ attains $P^{Monge}_{\mu,\nu}$. Conversely, assuming $\gamma_c= j(U)$, for some $U\in \mathcal{R}_a(\mu,\nu)$. Since $W$ is Polish, we can always find a Borel measurable map $F : W\to W$, whose $\mu-$equivalence class, as a $\mathcal{B}(W)^\mu-$measurable map, is $U$. In particular, $\gamma_c = (I_W\times F)_\star \mu$ ; $F$ meets the assumptions of a unique strong solution, for the associated equation. Finally, the symmetric counterpart, follows similarly.  \nqed

\end{document}

\end{document}